\documentclass[fleqn]{mat01}
\usepackage{times,mathtimy,amssymb,latexsym}
\begin{document}

\setcounter{page}{319}
\firstpage{319}

\font\zz=msam10 at 10pt
\def\Box{\mbox{\zz{\char'244}}}

%\font\zz=msam10 at 10pt
%\def\BBox{\mbox{\zz{\char'245}}}

\font\azz=tibi at 10.4pt
\def\Sa{\mbox{\azz{S}}}
\def\L{\mbox{\azz{L}}}

\newtheorem{theore}{Theorem}
\renewcommand\thetheore{\arabic{section}.\arabic{subsection}.\arabic{theore}}
\newtheorem{theor}[theore]{\bf Theorem}

\newtheorem{propo}[theore]{\rm PROPOSITION}
\newtheorem{rem}[theore]{Remark}
\newtheorem{coro}[theore]{\rm COROLLARY}
\newtheorem{definit}[theore]{\rm DEFINITION}

\renewcommand{\theequation}{\thesubsection.\arabic{equation}}

\title{Conductors and newforms for ${\hbox{\bfseries\itshape U}\hbox{\bf (1,1)}}$}

\markboth{Joshua Lansky and A~Raghuram}{Newforms for $U(1,1)$}

\author{JOSHUA LANSKY and A~RAGHURAM$^{*}$}

\address{Department of Mathematics, American University, Washington DC 20910, USA\\
\noindent $^{*}$Department of Mathematics, University of Iowa, 14
Maclean Hall, Iowa City, IA~52242, USA\\
\noindent E-mail: lansky@american.edu; araghura@math.uiowa.edu}

\volume{114}

\mon{November}

\parts{4}

\Date{MS received 18 May 2004; revised 23 June 2004}

\begin{abstract}
Let $F$ be a non-Archimedean local field whose residue characteristic is
odd. In this paper we develop a theory of newforms for $U(1,1)(F)$,
building on previous work on $SL_2(F)$. This theory is analogous to the
results of Casselman for $GL_2(F)$ and Jacquet, Piatetski-Shapiro, and
Shalika for $GL_n(F)$. To a representation $\pi$ of $U(1,1)(F)$, we
attach an integer $c(\pi)$ called the conductor of $\pi$, which depends
only on the $L$-packet $\Pi$ containing $\pi$. A newform is a vector in
$\pi$ which is essentially fixed by a congruence subgroup of level
$c(\pi)$. We show that our newforms are always test vectors for some
standard Whittaker functionals, and, in doing so, we give various
explicit formulae for newforms.
\end{abstract}

\keyword{Conductor; newforms; representations; $U(1,1)$.}

\maketitle

\section{Introduction}

To introduce the main theme of this paper we recall the following
theorem of Casselman \cite{casselman}. Let $F$ be a non-Archimedean
local field whose ring of integers is $\mathcal{O}_F$. Let
$\mathcal{P}_F$ be the maximal ideal of $\mathcal{O}_F$. Let $\psi_F$ be
a non-trivial additive character of $F$ which is normalized so that the
maximal fractional ideal on which it is trivial is $\mathcal{O}_F$.

\begin{theor}[\cite{casselman}] \label{thm:casselman}
Let $(\pi,V)$ be an irreducible admissible infinite-dimensional
representation of $GL_2(F)$. Let $\omega_{\pi}$ denote the central
character of $\pi$. Let
\begin{equation*}
\Gamma(m) = \left\{ \left( \begin{array}{cc}
a &b\\
c &d
\end{array} \right) \in GL_2 (\mathcal{O}_F) :
c \equiv 0 \pmod{\mathcal{P}_F^m} \right\}.
\end{equation*}
Let
\begin{equation*}
V_m = \left\{v \in V : \pi \left( \left( \begin{array}{cc}
a &b \\
c &d \end{array}\right) \right) v = \omega_{\pi}(d)v,\quad
\forall \left( \begin{array}{cc}
a &b \\
c &d \end{array}\right) \in \Gamma(m) \right\}.
\end{equation*}

$\left.\right.$\vspace{-2pc}
\begin{enumerate}
\renewcommand\labelenumi{\rm (\roman{enumi})}
\leftskip .2pc
\item There exists a non-negative integer $m$ such that $V_m \not=(0).$
If $c(\pi)$ denotes the least non-negative integer $m$ with this
property then the epsilon factor $\epsilon(s,\pi,\psi_F)$ of $\pi$ is up
to a constant multiple of the form $q^{-c(\pi)s}$. {\rm (}Here $q$ is the
cardinality of the residue field of $F$.{\rm )}

\item For all $m \geq c(\pi)$ we have $\dim (V_m) = m-c(\pi)+1$.
\end{enumerate}
\end{theor}

The assertion $\dim (V_{c(\pi)}) = 1$ is sometimes referred to as {\it
multiplicity one theorem for newforms} and the unique vector (up to
scalars) in $V_{c(\pi)}$ is called the {\it newform} for $\pi.$ This is
closely related to the classical Atkin--Lehner theory of newforms for
holomorphic cusp forms on the upper half plane \cite{casselman}. When
$c(\pi) = 0$ we have a spherical representation and the newform is
nothing but the spherical vector.

Newforms play an important role in the theory of automorphic forms. We
cite two examples to illustrate this. First, the zeta integral
corresponding to the newform is exactly the local $L$-factor associated
to $\pi$ (see \cite{jacquet-ps-shalika} for instance). In addition,
newforms frequently play the role of being `test vectors' for
interesting linear forms associated to $\pi$. For example, the newform
is a test vector for an appropriate Whittaker linear functional. In
showing this, explicit formulae for newforms are quite often needed. For
instance, if $\pi$ is a supercuspidal representation which is realized
in its Kirillov model then the newform is the characteristic function of
the unit group $\mathcal{O}_F^{\times}$. This observation is implicit in
Casselman \cite{casselman} and is explicitly stated and proved in
Shimizu \cite{schmidt}. Since the Whittaker functional on the Kirillov
model is given by evaluating functions at $1\in F^*$, we get in
particular that the functional is non-zero on the newform. In a related
vein \cite{dipendra} and \cite{gross-prasad} show that test vectors for
trilinear forms for $GL_2(F)$ are often built from newforms. (See
also a recent expository paper of Schmidt \cite{schmidt} where many of
these results are documented.)

In addition to Casselman's theory for $GL_2(F)$, newforms have been
studied for certain other classes of groups. Jacquet
{\it et~al}~\cite{jacquet-ps-shalika} have developed a theory of newforms for
{\it generic} representations of $GL_n(F)$. In this setting, there is no
satisfactory statement analogous to (ii) of the above theorem. However,
in his recent thesis, Mann \cite{12} obtained several results
on the growth of the dimensions of spaces of fixed vectors and has a
conjecture about this in general. For the group $GL_2(D)$, $D$ a
$p$-adic division algebra, Prasad and Raghuram \cite{dipendra-raghuram}
have proved an analogue of Casselman's theorem for irreducible principal
series representations and supercuspidal representations coming via
compact induction. In an unpublished work, Brooks Roberts has proved
part of (i) of the above \hbox{theorem} for representations of $GSp_4(F)$ whose
Langlands parameter is induced from a two-dimensional representation of
the Weil--Deligne group of $F.$ In a previous paper
\cite{josh-raghuram2}, we develop a theory of conductors and newforms
for $SL_2(F)$. In this paper we use the results of \cite{josh-raghuram2}
to carry out a similar program for the unramified quasi split unitary
group\break $U(1,1)$.

Let $\bar{G} = U(1,1)(F).$ Crucial to our study of newforms are certain
filtrations of maximal compact subgroups of $\bar{G}.$ Let $\bar{K} =
\bar{K}_0$ be the standard hyperspecial maximal compact subgroup of
$\bar{G}$. Let $\bar{K}' = \bar{K}_0' = \alpha^{-1} \bar{K}_0 \alpha$,
where $\alpha = \left(\begin{smallmatrix} \varpi_F & 0 \\ 0 & 1
\end{smallmatrix}\right).$ Then $\bar{K}_0$ and $\bar{K}_0'$ are, up to
conjugacy, the two maximal compact subgroups of $\bar{G}$. We define
filtrations of these maximal compact subgroups as follows. For $m$ an
integer $\geq 1,$ let
\begin{equation*}
\bar{K}_m = \left\{ \left( \begin{array}{cc}
a &b\\
c &d
\end{array} \right) \in \bar{K}_0\hbox{:}\ c \equiv 0
\pmod{\mathcal{P}_E^m} \right\} \quad {\rm and} \quad \bar{K}_m' =
\alpha^{-1}\bar{K}_m \alpha.
\end{equation*}

$\left.\right.$\vspace{-1.5pc}

Let $(\bar{\pi},V)$ be an irreducible admissible infinite-dimensional
representation of $\bar{G}.$ Let $\bar{Z}$ denote the center of
$\bar{G}$ and let $\omega_{\bar{\pi}}$ be the central character of
$\bar{\pi}$. Let $\bar{\eta}$ be any character of
$\mathcal{O}_E^{\times}$ such that $\bar{\eta} = \omega_{\bar{\pi}}$ on
the center. Let $c(\bar{\eta})$ denote the conductor of $\bar{\eta}.$
For any $m \geq c(\bar{\eta}),$ $\bar{\eta}$ gives a character of
$\bar{K}_m$ and also $\bar{K}_m'$ given by
$\bar{\eta} \left( \left( \begin{smallmatrix} a &b\\ c &d
\end{smallmatrix} \right) \right) = \bar{\eta}(d). $ We define for $m \geq
c(\bar{\eta})$,
\begin{equation*}
\bar{\pi}_{\bar{\eta}}^{\bar{K}_m} := \left\{v \in V\hbox{:}\ \bar{\pi}
\left( \left( \begin{array}{cc}
a &b \\
c &d \end{array} \right) \right)v = \bar{\eta}(d)v,\quad
\forall \left( \begin{array}{cc}
a &b \\
c &d \end{array} \right) \in \bar{K}_m \right\}.
\end{equation*}
The space $\bar{\pi}_{\bar{\eta}}^{\bar{K}_m'}$ is defined analogously.
We define the {\it $\bar{\eta}$-conductor} $c_{\bar{\eta}}(\bar{\pi})$
of $\bar{\pi}$ as
\setcounter{equation}{1}
\begin{equation}
c_{\bar{\eta}} (\bar{\pi}) = \min\{m \geq 0\hbox{:}\
\bar{\pi}_{\bar{\eta}}^{\bar{K}_m} \not= (0)
\quad {\rm or} \quad \bar{\pi}_{\bar{\eta}}^{\bar{K}_m'} \not= (0)\}.
\end{equation}

We define the {\it conductor $c(\bar{\pi})$} of $\bar{\pi}$ by
\begin{equation}
c(\bar{\pi})  =  \min \{c_{\bar{\eta}}(\bar{\pi}) : \bar{\eta}\},
\end{equation}
where $\bar{\eta}$ runs over characters of $\mathcal{O}_E^{\times}$
which restrict to the central character $\omega_{\bar{\pi}}$ on
$\bar{Z}$.

We deal with the following basic issues in this paper.
\begin{enumerate}
\renewcommand\labelenumi{(\roman{enumi})}
\leftskip .3pc
\item Given an irreducible representation $\bar{\pi}$, we determine its
conductor $c(\bar{\pi}).$ A very easy consequence (almost built into the
definition) is that the conductor depends only on the $L$-packet
containing $\bar{\pi}$.

\item We identify the conductor with other invariants associated to the
representation. For instance, for $SL_2(F)$ we have shown
\cite{josh-raghuram2} that the conductor of a representation is same as
the conductor of a minimal representation of $GL_2(F)$ determining its
$L$-packet. We prove a similar result for $U(1,1)(F)$ in this paper. See
\S\ref{sec:comparison} and \S\ref{sec:comparison-u11}.

\item We determine the growth of the space ${\rm
dim}(V^{\bar{K}_m}_{\bar{\eta}})$ as a function of $m$. This question is
analogous to (ii) of Casselman's theorem quoted above. Computing such
dimensions is of importance in `local level raising' issues. See
\cite{12}.

\item We address the question of whether there is a multiplicity one
result for newforms. It turns out that quite often ${\rm
dim}(V^{\bar{K}_{c(\bar{\pi})}}_{\bar{\eta}}) = 1$, but this fails in
general (for principal series representations of a certain kind). In
these exceptional cases the dimension of the space of newforms is two,
but a canonical quotient of this two-dimensional space has dimension one
(see~\S\ref{sec:multiplicity-one}).

\item Are the newforms test vectors for Whittaker functionals? This is
important in global issues related to newforms. We are grateful to
Benedict Gross for suggesting this question to us. It turns out that our
newforms are always test vectors for Whittaker functionals. In the
proofs we often need explicit formulae for newforms in various models
for the representations. These formulas are interesting for their own
sake. For example, if $(\pi,V)$ is a ramified supercuspidal
representation of $U(1,1)(F)$, then the newform can be taken as the
characteristic function of $(\mathcal{O}_F^{\times})^2$ where $V$ is
regarded as a subspace of the Kirillov model of a canonically associated
minimal representation of $GL_2(F)$ (cf. \cite{schmidt}).
\end{enumerate}

We set up notation in \S\ref{sec:notation} following that used in
\cite{josh-raghuram2}. We then briefly review the structure of
$L$-packets for $SL_2$ and $U(1,1)$ in \S\ref{sec:packets-sl2-u11}. As
this paper depends crucially on our previous paper~\cite{josh-raghuram2}
on $SL_2$, we summarize the results of~\cite{josh-raghuram2} in
\S\ref{sec:newforms}. The heart of this paper is \S\ref{sec:unitary}. In
\S\ref{sec:defns-u11} we define the notion of conductor and then make
some easy but technically important remarks on spaces of fixed vectors.
The next two subsections deal respectively with sub-quotients of
principal series representations and supercuspidal representations. In
\cite{josh-raghuram2}, we use Kutzko's construction of supercuspidal
representations of $GL_2(F)$ to obtain results for supercuspidals of
$SL_2(F)$. In this paper, we use these results, in turn, to obtain
information for $U(1,1)(F)$. In general, we will often reduce the proofs
of statements concerning $U(1,1)(F)$ to those of the corresponding
$SL_2(F)$ statements. In particular, we exploit the fact that $SL_2(F)$
is the derived group of $U(1,1)(F)$ and that $U(1)(F)SL_2(F)$ has
index two in $U(1,1)(F).$ In this way we avoid directly dealing with
$K$-types and other intrinsic details for ${\rm U}(1,1)(F)$ as much of
the work has been done for $SL_2(F)$ in~\cite{josh-raghuram2}. Finally,
in~\S\ref{sec:multiplicity-one}, we briefly discuss a multiplicity one
result for\break newforms.

We mention some further directions that arise naturally from this work.
To begin with, it would be interesting to see how our theory of newforms
and conductors bears upon known results about local factors for
$U(1,1)(F)$. In particular, are our conductors the same as (or closely
related to) the analytic conductors appearing in epsilon factors? Also,
is a zeta-integral corresponding to a newform of a representation equal
to a local $L$-factor for the representation?

\section{Preliminaries}\label{sec:prelims}

\subsection{\it Notation}\label{sec:notation}

If $L$ is any non-Archimedean local field let $\mathcal{O}_L$ be its
ring of integers and let $\mathcal{P}_L$ be the maximal ideal of
$\mathcal{O}_L.$ Let $\varpi_L$ be a uniformizer for $L$, i.e.,
$\mathcal{P}_L = \varpi_L\mathcal{O}_L.$ Let $k_L =
\mathcal{O}_L/\mathcal{P}_L$ be the residue field of $L.$ Let $p$ be the
characteristic of $k_L$ and let the cardinality of $k_L$ be $q_L$ which
is a power of $p.$ Let $\epsilon_L$ be an element of $\mathcal{O}_L^* -
\mathcal{O}_L^{*2}$.

If $n$ is a positive integer, let $U_L^n$ denote the $n$th filtration
subgroup $1+\mathcal{P}^n_L$ of $\mathcal{O}_L^{\times}$, and define
$U_L^0 = \mathcal{O}_L^{\times}$. Let $\mathfrak{v}_L$ denote the
additive valuation on $L^*$ which takes the value $1$ on $\varpi_L.$ We
let $|\cdot |_L$ denote the normalized multiplicative valuation given by
$|x|_L = q^{-\mathfrak{v}_L(x)}.$ If $\chi$ is a character of $L^*$ we
define the conductor $c(\chi)$ to be the smallest non-negative integer
$n$ such that $\chi$ is trivial on $U_L^n$. Let $\psi_L$ be a
non-trivial additive character of $L$ which is assumed to be trivial on
$\mathcal{O}_L$ and non-trivial on $\mathcal{P}_L^{-1}.$ For any $a \in
L$ the character given by sending $x$ to $\psi_L(ax)$ will be denoted as
$\psi_{L,a}$ or simply by $\psi_a.$ (In all the above notations we may
omit the subscript $L$ if there is only one field in the context.)

In the following, $F$ will be a fixed non-Archimedean local field whose
residue characteristic is odd and $E$ will be used to denote a quadratic
extension of $F$. We denote by $\omega_{E/F}$ the quadratic character of
$F^*$ associated to $E/F$ by local class field theory. Recall that the
kernel of $\omega_{E/F}$ is $N_{E/F}(E^*)$, the norms from $E^*.$ We
will require the units $\epsilon_F$ and $\epsilon_E$ to be compatible in
the sense that
\begin{equation*}
\epsilon_F = N_{E/F}(\epsilon_E).
\end{equation*}

We let $\widetilde{G}$ denote the group $GL_2(F)$. Let $\widetilde{B} =
\widetilde{T}N$ be the standard Borel\vspace{.1pc} subgroup of upper triangular
matrices in $\widetilde{G}$ with Levi subgroup $\widetilde{T}$ and
unipotent radical $N$. Let $\widetilde{Z}$ be the center of
$\widetilde{G}$. Let $G = SL_2(F).$ Let $B = TN$ be the standard Borel
subgroup of upper triangular matrices in $G$ with Levi subgroup $T$ and
unipotent radical $N$. We set $K=SL_2(\mathcal{O}_F)$ and
$\widetilde{K}=GL_2(\mathcal{O}_F)$ and denote by $I$ and
$\widetilde{I}$ respectively the standard Iwahori subgroups of $G$ and
$\widetilde{G}$.

Suppose that $E/F$ is unramified, and let $s$ denote the non-trivial
element of $\text{Gal}(E/F)$. We denote by $\bar{G}$ the group $U(1,1)$,
i.e., the group of all $g\in GL_2(E)$ such that
\begin{equation*}
{^{s}\!g} \left(\begin{array}{cc}
0 &1\\
-1 &0 \end{array}\right)\, {^{t}\!g} = \left(\begin{array}{cc}
0 &1\\
-1 &0 \end{array}\right).
\end{equation*}
We let $\bar{B}$ be the standard upper triangular Borel subgroup of
$\bar{G}$ with diagonal Levi subgroup $\bar{T}$ and unipotent radical
$N$. We note that the elements of $\bar{T}$ are of the form
$\left(\begin{smallmatrix} t &0\\ 0 &^{s}t^{-1}
\end{smallmatrix}\right)$ for $t\in E^*$, and those of $\bar{B}$ are of
the form $\left(\begin{smallmatrix} t & ta\\ 0 &^{s}t^{-1}
\end{smallmatrix}\right)$ with $t \in E^*$ and $a \in F.$ We let
$\bar{Z}$ be the center of $\bar{G}$. Then $\bar{Z}\cong E^1$, where
$E^1 = {\rm ker}(N_{E/F})$ is the subgroup of norm one elements of
$E^*.$ Denote by $\bar{I}$ the standard Iwahori subgroup of $\bar{G}$
and by $\bar{K}$ the standard hyperspecial maximal compact subgroup of
$\bar{G}$.

The following filtrations of maximal compact subgroups of $G$ will be
important in our study of newforms. Let $K_{-1} = G$ and $K_0 = K$. Let
$K' = K_0' = \alpha^{-1} K_0 \alpha$, where $\alpha = \left(
\begin{smallmatrix} \varpi_F &0\\ 0 &1 \end{smallmatrix} \right).$ Then
$K_0$ and $K_0'$ are, up to conjugacy, the two maximal compact subgroups
of $G$. For $m$ an integer $\geq 1,$
\begin{align*}
K_m &= \left\{ \left(\begin{array}{cc}
a &b \\
c &d \end{array}\right) \in K_0 :
c \equiv 0 \pmod{\mathcal{P}_F^m} \right\},\\
K_m' &= \alpha^{-1}K_m \alpha.
\end{align*}
We note that for $m \geq 1$ the following inclusions hold up to
conjugacy within $G$:
\setcounter{equation}{0}
\begin{equation}
K_{m+1}' \subset K_m \subset K_{m-1}'.
\end{equation}
Analogous results hold for the following filtration groups of $\bar{G}$:
\begin{align*}
\bar{K}_{-1} &= \bar{G},\\
\bar{K}_0 &= \bar{K},\\
\bar{K}_m &= \left\{ \left(\begin{array}{cc}
a &b\\
c &d \end{array}\right) \in \bar{K}_0 :
c \equiv 0 \pmod{\mathcal{P}_E^m} \right\},\\
\bar{K}_m' &= \alpha^{-1}\bar{K}_m \alpha.
\end{align*}
We note that the filtration subgroups for $G$ and $\bar{G}$ are related
by
\begin{equation}\label{eq:filtration}
\bar{K}_m = K_m \bar{T}_0,
\end{equation}
where $\bar{T}_0 = \bar{T} \cap \bar{K}_0.$

In addition to $\alpha$, we will also make frequent use of the matrices
$\beta := \left(\begin{smallmatrix}
1 &0\\ 0 &\varpi_F
\end{smallmatrix}\right),$  $\gamma := \left(\begin{smallmatrix}
\epsilon_F &0\\ 0 &1\end{smallmatrix}\right)$ and
$\theta := \left(\begin{smallmatrix}
\epsilon_E &0\\ 0 &^{s}\epsilon_E^{-1}\end{smallmatrix}\right).$

For any subsets $A,B,C,D \subset F$ we let
\begin{equation*}
\left[ \begin{array}{cc}
A &B\\
C &D
\end{array}\right] = \left\{ \left( \begin{array}{cc}
a &b\\
c &d
\end{array}\right) : a \in A, b \in B, c \in C, d \in D \right\}.
\end{equation*}
We denote $\left[ \begin{smallmatrix} 1 &\mathcal{P}^j\\ 0 &1
\end{smallmatrix}\right]$ by $N(\mathcal{P}^j)$ or simply by $N(j).$ We
let $\overline{N}$ denote the lower triangular unipotent subgroup of $G$
and a similar meaning is given to $\overline{N}(\mathcal{P}^j)$ and
$\overline{N}(j).$

If $\mathcal{H}$ is a closed subgroup of a locally compact group
$\mathcal{G}$ and if $\sigma$ is an admissible representation of
$\mathcal{H}$ then ${\rm Ind}_{\mathcal{H}}^{\mathcal{G}}(\sigma)$
denotes {\it normalized} induction, and ${\rm
ind}_{\mathcal{H}}^{\mathcal{G}}(\sigma)$ denotes the subrepresentation
of ${\rm Ind}_{\mathcal{H}}^{\mathcal{G}}(\sigma)$ consisting of those
functions whose support is compact mod $\mathcal{H}.$ The symbol ${\bf
1}$ will denote the trivial representation of the group in context.

For any real number $\zeta$ we let $\lceil \zeta \rceil$ denote the
least integer greater than or equal to $\zeta$ and we let $\lfloor \zeta
\rfloor = - \lceil -\zeta \rceil.$

\subsection{\it $L$-packets for $SL_2(F)$ and $U(1,1)$}\label{sec:packets-sl2-u11}

In this section we collect statements about the structure of $L$-packets
for $G=SL_2(F)$ and $\bar{G} = U(1,1)$. All the assertions made here are
well-known and can be read off from a combination of Labesse and Langlands
\cite{labesse-langlands}, Gelbart and Knapp \cite{gelbart-knapp} and
Rogawski \cite{rogawski}.

If $\widetilde{\pi}$ is an irreducible admissible representation of
$\widetilde{G}$ then its restriction to $G$ is a multiplicity free
finite direct sum of irreducible admissible representations of $G$ which
we often write as
\begin{equation*}
{\rm Res}_{SL_2(F)} \widetilde{\pi} = \pi_1 \oplus \cdots \oplus \pi_r.
\end{equation*}

On the other hand, if $\pi$ is any irreducible admissible representation
of $G$ then there exists an irreducible admissible representation
$\widetilde{\pi}$ of $\widetilde{G}$ whose restriction to $G$ contains
$\pi.$

Note that $\widetilde{G}$ acts on the space of all equivalence classes
of irreducible admissible representations of $G$ and an {\it $L$-packet}
for $G$ is simply an orbit under this action. It turns out that, with
the notation as above, the $L$-packets are precisely the sets
$\{\pi_1,\ldots,\pi_r\}$ appearing in the restrictions of irreducible
representations $\widetilde{\pi}$ of $\widetilde{G}.$

We now give some general statements concerning the $L$-packets for
$\bar{G} = U(1,1)$. The adjoint group of $U(1,1)$ is $PGL_2$, and hence
$PGL_2(F)$ and $\widetilde{G}$ act via automorphisms on $\bar{G}$, hence
act on the set of equivalence classes of irreducible representations of
$\bar{G}$. Rogawski~(\cite{rogawski}, \S11.1) defines an $L$-packet for
$\bar{G}$ to be an orbit under this action. If $\bar{\pi}$ is an element
of a non-trivial $L$-packet, then the other element of the $L$-packet is
$^\alpha \bar{\pi}$.

If $\bar{\Pi}$ is an $L$-packet for $\bar{G}$, then the set of
irreducible components of the restrictions of elements of $\bar{\Pi}$ to
$G$ is an $L$-packet $\Pi$ for $G$. The direct sum
$\bigoplus_{\pi\in\Pi}\pi$ is therefore the restriction of an
irreducible admissible representation $\widetilde{\pi}$ of
$\widetilde{G}.$ This $\widetilde{\pi}$ is unique up to twisting by a
character. In practice, we will choose a convenient $\widetilde{\pi}.$
Since $\bigoplus_{\pi\in\Pi}\pi = \text{Res}_G
\left(\bigoplus_{\bar{\pi}\in\bar{\Pi}}\bar{\pi}\right)$, we obtain an
action of $\widetilde{G}$ on
$\bigoplus_{\bar{\pi}\in\bar{\Pi}}\bar{\pi}$ via the represen-\break
tation $\widetilde{\pi}$.

\section{Newforms for $\Sa\L_{\bf 2}$}\label{sec:newforms}

This section collects our results \cite{josh-raghuram2} on conductors
and newforms for $SL_2(F)$. All these results, along with their complete
proofs, can be found in \cite{josh-raghuram2}.

\subsection{\it Definitions}\label{sec:defn-sl2}

We now give our definition of the conductor of a representation of $G$.
The basic filtration subgroups of $G$ considered in this paper are $K_0
= K = SL_2(\mathcal{O}_F)$ and for $m \geq 1$,
\begin{equation*}
K_m = \left\{ \left(\begin{array}{cc}
a &b \\
c &d \end{array}\right) \in SL_2(\mathcal{O}_F) :
c \equiv 0 \pmod{\mathcal{P}_F^m} \right\}.
\end{equation*}
For all $m \geq 0$ we let $K_m' = \alpha^{-1}K_m\alpha.$

Let $(\pi,V)$ be an irreducible admissible infinite-dimensional
representation of $G.$ Let $\omega_{\pi}$ be the character of $\{\pm
1\}$ such that $\pi\left(\left(\begin{smallmatrix} -1 &0\\ 0 &-1
\end{smallmatrix}\right)\right) = \omega_{\pi}(-1)1_V.$

We let $\eta$ be any character of $\mathcal{O}_F^\times$ such that
$\eta(-1) = \omega_{\pi}(-1).$ Let $c(\eta)$ denote the conductor of
$\eta.$ For any $m \geq c(\eta),$ $\eta$ gives a character of $K_m$ and
also $K_m'$ given by
\begin{equation*}
\eta\left(\left(\begin{matrix}
a &b\\
c &d \end{matrix}\right)\right)
= \eta(d).
\end{equation*}

We define
\begin{equation*}
\pi_{\eta}^{K_m} := \left\{v \in V : \pi \left( \left( \begin{array}{cc}
a &b\\
c &d \end{array} \right) \right) v = \eta(d)v, \quad
\forall \left( \begin{array}{cc}
a &b\\
c &d \end{array} \right) \in K_m \right\}.
\end{equation*}
The spaces $\pi_{\eta}^{K_m'}$ are defined analogously. We note that
$\pi_{\eta}^{K_m} = \pi_{\eta}^{K'_m} = (0)$ for $m < c(\eta)$.

We define the {\it $\eta$-conductor} $c_{\eta}(\pi)$ of $\pi$ as
\setcounter{equation}{0}
\begin{equation}
c_{\eta}(\pi) = \min\{m \geq 0 : \pi_{\eta}^{K_m} \not = (0) \quad {\rm
or} \quad \pi_{\eta}^{K_m'} \not= (0)\},
\end{equation}

We define the {\it conductor $c(\pi)$} of $\pi$ by
\begin{equation}
c(\pi)  =  \min \{c_{\eta}(\pi) : \eta \},
\end{equation}
where $\eta$ runs over characters of $\mathcal{O}_F^\times$ such that
$\eta(-1) = \omega_{\pi}(-1).$ If $\eta$ is such that $c_\eta(\pi) =
c(\pi)$ and $\pi_{\eta}^{K_{c(\pi)}}\neq (0)$
(resp.~$\pi_{\eta}^{K'_{c(\pi)}}\neq (0)$), then we call
$\pi_\eta^{K_{c(\pi)}}$ (resp.~$\pi_\eta^{K'_{c(\pi)}}$) a \textit{space
of newforms} of $\pi$. In this case, we refer to a non-zero element of
$\pi_{\eta}^{K_{c(\pi)}}$ or $\pi_{\eta}^{K'_{c(\pi)}}$ as a
\textit{newform} of $\pi$.

\subsection{\it Principal series representations}\label{sec:principal-sl2}

Let $\chi$ be a character of $F^*.$ Then $\chi$ gives a character of
$B$ via the formula $\chi\left(\left(\begin{smallmatrix}
a &b\\
 &a^{-1} \end{smallmatrix}\right)\right) = \chi(a).$
Let $\pi(\chi)$ stand for the (unitarily) induced representation ${\rm
Ind}_B^G(\chi).$ It is well-known that $\pi(\chi)$ is reducible if and
only if $\chi$ is either $|\cdot |_F^{\pm}$ or if $\chi$ is a quadratic
character.

There is an essential difference between the two kinds of
reducibilities. If $\chi = |\cdot |_F^{\pm}$, then $\pi(\chi)$ is the
restriction to $G$ of a reducible principal series representation of
$\widetilde{G}$. Hence $\pi(\chi)$ will have two representations in its
Jordan--H\"older series, namely the trivial representation and the
Steinberg representation which we will denote by ${\rm St}_G.$

If $\chi$ is a quadratic character, then $\pi(\chi)$ is the restriction
to $G$ of an irreducible principal series representation of
$\widetilde{G}$ and breaks up as a direct sum of two irreducible
representations, which constitute an $L$-packet of $G$. If $\chi =
\omega_{E/F}$ we denote $\pi(\chi)$ by $\pi_E$ and let $\pi_E \simeq
\pi_E^1 \oplus \pi_E^2.$ We denote the $L$-packet by $\xi_E = \{\pi_E^1,
\pi_E^2 \}.$

As mentioned in the introduction, one of the applications of newforms we
have in mind is that they are test vectors for Whittaker functionals.
For principal series representations and in fact all their sub-quotients
we consider the following $\psi$-Whittaker functional (see
\cite{schmidt}). For any function $f$ in a principal series
representation $\pi(\chi)$ we define
\setcounter{equation}{0}
\begin{equation}\label{eqn:whittaker-ps}
\Lambda_{\psi}f := \underset{r \to \infty}{\lim}
\int_{\mathcal{P}_F^{-r}} f \left(
\left(\begin{array}{cc} 0 &-1\\ 1 &0 \end{array}\right)
\left(\begin{array}{cc} 1 &x\\ 0 &1 \end{array}\right) \right)
\overline{\psi}(x)\, {\rm d}x,
\end{equation}
where the Haar measure is normalized such that ${\rm vol}(\mathcal{O}) =
1.$

\begin{propo} {\rm (Unramified principal series representations).}\label{prop:ps-sl2}
$\left.\right.$\vspace{.5pc}

\noindent Let $\chi$ be an unramified character of $F^*$ and let
$\pi(\chi)$ be the corresponding principal series representation of $G.$
We have
\begin{equation*}
\dim (\pi(\chi)^{K_m}) = \left\{ \begin{array}{ll}
1, &\mbox{ if $m=0,$}\\[.2pc]
2m, &\mbox{ if $m \geq 1.$} \end{array} \right.
\end{equation*}
\end{propo}

\begin{coro} {\rm (Test vectors for unramified principal series representations)}\label{cor:unramified-sl2-test}
$\left.\right.$\vspace{.5pc}

\noindent For an unramified character $\chi$ of $F^*$ such that $\chi
\not= |\cdot |_F^{-1}.$ Let $f_{\rm new}$ be any non-zero $K$-fixed
vector. Then we have
\begin{equation*}
\Lambda_{\psi}f_{\rm new} = L(1,\chi)^{-1} \not= 0,
\end{equation*}
where $L(s,\chi)$ is the standard local abelian $L$-factor associated to
$\chi.$\vspace{.2pc}
\end{coro}

\begin{propo} {\rm (Steinberg representation)}\label{prop:stienberg-sl2}
$\left.\right.$\vspace{.5pc}

\noindent If \ ${\rm St}_G$ is the Steinberg representation of $G,$ then
the dimension of the space of fixed vectors under $K_m$ is given by
\begin{equation*}
{\rm dim}({\rm St}_G)^{K_m} =  \left\{ \begin{array}{ll}
0, &\mbox{ if $m=0,$}\\[.2pc]
2m-1, &\mbox{ if $m \geq 1.$} \end{array} \right.
\end{equation*}
\end{propo}\vspace{.2pc}

\begin{coro} {\rm (Test vectors for the Steinberg representation)}\label{cor:steinberg-test}
$\left.\right.$\vspace{.5pc}

\noindent Let the Steinberg representation ${\rm St}_G$ be realized as
the unique irreducible subrepresentation of $\pi(|\cdot |).$ The
$\psi$-Whittaker functional $\Lambda_{\psi}$ is non-zero on the space of
newforms $({\rm St}_G)_{\rm new} = {\rm St}_G^{K_1}$.
\end{coro}\vspace{.2pc}

\begin{propo} {\rm (Ramified principal series representations)}\label{prop:sl2-ps-ramified}
$\left.\right.$\vspace{.5pc}

\noindent Let $\chi$ be a ramified character of $F^*.$ Let $\pi =
\pi(\chi)$ be the corresponding principal series representation of $G.$
Let $c(\chi)$ denote the conductor of $\chi.$

\begin{enumerate}
\renewcommand\labelenumi{\rm (\roman{enumi})}
\leftskip .3pc
\item We have $c(\pi) = c(\chi)$ and further $c_{\eta}(\pi) = c(\pi)$
only for those characters $\eta$ such that $\eta = \chi^{\pm}$ on the
group of units $\mathcal{O}^{\times}.$

\item If $\chi^2|_{(\mathcal{O}^{\times})^2} \not= {\bf 1}$ and $\eta =
\chi |_{\mathcal{O}^{\times}}$ then
\begin{equation*}
\hskip -1.25pc\dim (\pi(\chi)_{\eta}^{K_m}) = \left\{ \begin{array}{ll}
0, &\mbox{ if $m < c(\chi),$}\\[.2pc]
1, &\mbox{ if $m = c(\chi),$}\\[.2pc]
2(m-c(\chi))+1, &\mbox{ if $m > c(\chi).$} \end{array} \right.
\end{equation*}

\item If $\chi^2|_{(\mathcal{O}^{\times})^2} = {\bf 1}$ and $\eta = \chi
|_{\mathcal{O}^{\times}}$ then
\begin{equation*}
\hskip -1.25pc\dim (\pi(\chi)_{\eta}^{K_m}) = \left\{ \begin{array}{ll}
0, &\mbox{ if $m = 0,$}\\[.2pc]
2m, &\mbox{ if $m \geq 1 = c(\chi)$.}
\end{array} \right.
\end{equation*}
\end{enumerate}
\end{propo}

\begin{coro} {\rm (Test vectors for ramified principal series representations)}
\label{cor:test-ram-ps}$\left.\right.$\vspace{.5pc}

\noindent Let $\chi$ be a ramified character of $F^*.$ Let $\pi =
\pi(\chi)$ be the corresponding principal series representation of $G.$
Assume that $\pi$ is irreducible. Let $m = c(\chi) \geq 1$ denote the
conductor of $\chi.$ The space of newforms $\pi(\chi)_{\rm new} =
\pi(\chi)^{K_{c(\chi)}}_{\chi}$ is one-dimensional and the Whittaker
functional $\Lambda_{\psi}$ is non-zero on this space of newforms.\vspace{.2pc}
\end{coro}

\begin{propo} {\rm (Ramified principal series $L$-packets)}\label{prop:ramified-sl2}
$\left.\right.$\vspace{.5pc}

\noindent Let $E/F$ be a quadratic ramified extension. Let $\xi_E =
\{\pi_E^1, \pi_E^2 \}$ be the corresponding $L$-packet. Then we have for
$\eta = \omega_{E/F},$\pagebreak
\begin{equation*}
\dim ((\pi_E^1)^{\eta}_m) = \dim ((\pi_E^2)^{\eta}_m) = \left\{ \begin{array}{ll}
0, &\mbox{ if $m = 0,$}\\[.2pc]
m, &\mbox{ if $m \geq 1.$} \end{array} \right.
\end{equation*}
\end{propo}

\begin{coro} {\rm (Test vectors for ramified principal series $L$-packets)}
\label{cor:test-ramified-l-packets-sl2}$\left.\right.$\vspace{.5pc}

\noindent Let $E/F$ be a ramified quadratic extension and let $\xi_E =
\{\pi_E^1, \pi_E^2 \}$ be the corresponding $L$-packet. Then one and
only one of the two representations in the packet is $\psi$-generic{\rm ,}
say{\rm ,} $\pi_E^1$. Then $\pi_E^2$ is $\psi_{\epsilon}$-generic. The
Whittaker functional $\Lambda_{\psi}$ is non-zero on the one dimensional
space of newforms $(\pi_E^1)_{\rm new} =
(\pi_E^1)_{\omega_{E/F}}^{K_1}.$ Any $\psi_{\epsilon}$-Whittaker
functional is non-zero on the one-dimensional space of newforms for
$\pi_E^2.$
\end{coro}

\begin{propo} {\rm (Unramified principal series $L$-packet)}
\label{prop:unramfified-sl2}$\left.\right.$\vspace{.5pc}

\noindent Let $E/F$ be the quadratic unramified extension. Let $\xi_E =
\{\pi_E^1, \pi_E^2\}$ be the corresponding $L$-packet.  Exactly one of
the two representations{\rm ,} say $\pi_E^1${\rm ,} has a non-zero
vector fixed by $K_0$. Then the dimensions of the space of fixed vectors
under $K_m$ and $K_m'$ for the two representations are as follows{\rm :}
\begin{enumerate}
\renewcommand\labelenumi{\rm (\roman{enumi})}
\leftskip .3pc
\item $\dim ((\pi_E^1)^{K_0}) = 1 = \dim ((\pi_E^2)^{K_0'}).$\vspace{.3pc}

\item $\dim ((\pi_E^1)^{K_0'}) = 0 = \dim ((\pi_E^2)^{K_0}).$\vspace{.3pc}

\item For $r \geq 1${\rm ,}
\begin{equation*}
\hskip -1.25pc\dim ((\pi_E^1)^{K_r}) = 2 \left\lfloor \frac{r}{2}
\right\rfloor + 1 = \dim ((\pi_E^2)^{K_r'}).
\end{equation*}
\item For $r \geq 1,$
\begin{equation*}
\hskip -1.25pc\dim ((\pi_E^1)^{K_r'}) = 2 \left\lfloor \frac{r-1}{2}
\right\rfloor + 1 = \dim ((\pi_E^2)^{K_r}).
\end{equation*}
\end{enumerate}
\end{propo}

\begin{coro} {\rm (Test vectors for unramified principal series $L$-packet)}
\label{cor:test-unramified-lpacket-sl2}$\left.\right.$\vspace{.5pc}

\noindent Let $E/F$ be the unramified quadratic extension{\rm ,} and let
$\xi_E = \{\pi_E^1, \pi_E^2 \}$ be the corresponding $L$-packet. Then
one and only one of the two representations in the packet is
$\psi$-generic{\rm ,} namely $\pi_E^1$ {\rm (}using the notation of
Proposition~{\rm \ref{prop:unramfified-sl2})}. Moreover, a
$\psi$-Whittaker functional is non-zero on the $K_0$-fixed vector in
$\pi_E^1$. The representation $\pi_E^2$ is not $\psi'$-generic for any
$\psi'$ of conductor $\mathcal{O}_F$. It is $\psi_{\varpi_F}$-generic
and any $\psi_{\varpi_F}$-Whittaker functional is non-zero on the unique
{\rm (}up to scalars{\rm )} $K_0'$-fixed vector in $\pi_E^2.$
\end{coro}

\subsection{\it Supercuspidal representations}\label{sec:sc-sl2}

We now consider supercuspidal representations of $G = SL_2(F).$ For this
we need some preliminaries on how they are constructed. We use Kutzko's
construction \cite{kutzko2,kutzko3} of supercuspidal representations for
$\widetilde{G}$ and then Moy and Sally \cite{moy-sally} or Kutzko and
Sally \cite{kutzko-sally} to obtain information on the supercuspidal
representations ($L$-packets) for $G.$

We begin by briefly recalling Kutzko's construction of supercuspidal
representations of $\widetilde{G}$ via compact induction from very
cuspidal representations of maximal open compact-mod-center subgroups.

For $l \geq 1$, let
\begin{equation*}
\widetilde{K}(l) = 1 + \mathcal{P}^lM_{2\times 2}(\mathcal{O})
\end{equation*}
be the principal congruence subgroup of $\widetilde{K}$ of level $l$.
Let $\widetilde{K}(0) = \widetilde{K}.$ Let $\widetilde{I}$ be the
standard Iwahori subgroup consisting of all elements in $\widetilde{K}$
that are upper triangular modulo $\mathcal{P}.$ For $l \geq 1$, let
\begin{equation*}
\widetilde{I}(l) = \left[\begin{array}{cc}
1 + \mathcal{P}^l &\mathcal{P}^l\\[.2pc]
\mathcal{P}^{l+1} &1 + \mathcal{P}^l \end{array}\right],
\end{equation*}
and let $\widetilde{I}(0) = \widetilde{I}.$ We will let $\widetilde{H}$
(resp. $\widetilde{J}$) denote either $Z\widetilde{K}$ (resp.
$\widetilde{K}$) or $N_{\widetilde{G}}\widetilde{I}$ (resp.
$\widetilde{I}$). Here $N_{\widetilde{G}}\widetilde{I}$ is the
normalizer in $\widetilde{G}$ of $\widetilde{I}.$ In either case we let
$\widetilde{J}(l)$ denote the corresponding filtration subgroup.

\setcounter{theore}{0}
\begin{definit}{\rm \cite{kutzko3,kutzko4}}\label{defn:very-cuspidal}
$\left.\right.$\vspace{.5pc}

\noindent {\rm An irreducible (and necessarily finite-dimensional)
representation $(\widetilde{\sigma},W)$ of $\widetilde{H}$ is called a
very cuspidal representation of level $l \geq 1$ if
\begin{enumerate}
\renewcommand\labelenumi{(\roman{enumi})}
\leftskip .1pc
\item $\widetilde{J}(l)$ is contained in the kernel of
$\widetilde{\sigma}.$

\item $\widetilde{\sigma}$ does not contain the trivial character
of $N(\mathcal{P}^{l-1}).$
\end{enumerate}}\vspace{-.5pc}
\end{definit}

We say that an irreducible admissible representation $\widetilde{\pi}$
of $\widetilde{G}$ is minimal if for every character $\chi$ of $F^*$
we have $c(\widetilde{\pi}) \leq c(\widetilde{\pi} \otimes \chi).$

\begin{theor}[\cite{kutzko3,kutzko4}] \label{thm:sc-kutzko}
There exists a bijective correspondence given by compact induction
$\widetilde{\sigma} \mapsto {\rm
ind}_{\widetilde{H}}^{\widetilde{G}}(\widetilde{\sigma})$ from very
cuspidal representations $\widetilde{\sigma}$ of either maximal open
compact-mod-center subgroup $\widetilde{H}$ and irreducible minimal
supercuspidal representations of $\widetilde{G}.$ Moreover{\rm ,} every
irreducible minimal supercuspidal representation of conductor $2l$
{\rm (}resp. $2l+1${\rm )} comes from a very cuspidal representation of
$Z\widetilde{K}$ {\rm (}resp. $N_{\widetilde{G}}\widetilde{I}${\rm )} of\break
level $l$.
\end{theor}

Following Kutzko we use the terminology that a supercuspidal
representation of $\widetilde{G}$ is said to be {\it unramified} if it
comes via compact induction from a representation of $Z\widetilde{K}$
and {\it ramified} if it comes via compact induction from a
representation of $N_{\widetilde{G}}\widetilde{I}$. We now take up both
types of supercuspidal representations and briefly review how they break
up on restriction to $G$. We refer the reader to \cite{kutzko-sally} and
\cite{moy-sally} for this.

We begin with the unramified case. Let $\widetilde{\sigma}$ be an
irreducible very cuspidal representation of $Z\widetilde{K}$ of level
$l$ ($\geq 1$). Let $\widetilde{\pi}$ be the corresponding supercuspidal
representation of $\widetilde{G}.$ Let $\sigma = {\rm
Res}_K(\widetilde{\sigma}).$ Then we have
\begin{equation*}
{\rm Res}_G(\widetilde{\pi}) = {\rm ind}_K^G(\sigma) \ \oplus \
^{\alpha}\!({\rm ind}_K^G(\sigma)),
\end{equation*}
where $\alpha = \left(\begin{smallmatrix} \varpi_F &0\\ 0 &1
\end{smallmatrix}\right).$

If $l \geq 2,$ or if $l=1$ and $\sigma$ is irreducible, then $\pi =
\pi(\sigma) = {\rm ind}_K^G(\sigma)$ is irreducible, hence so is $\pi' =
^\alpha\!\pi$. We thus have an unramified supercuspidal $L$-packet
$\{\pi,\pi'\}.$

If $l=1$ and $\sigma$ is reducible, then $\widetilde{\sigma}$ comes from
the unique (up to twists) cuspidal representation of $GL_2({\mathbb
F}_q)$ whose restriction to $SL_2({\mathbb F}_q)$ is reducible and
breaks up into the direct sum of the two cuspidal representations of
$SL_2({\mathbb F}_q)$ of dimension $(q-1)/2.$ Correspondingly, we have
$\sigma = \sigma_1 \oplus \sigma_2$, and if we let $\pi_i = {\rm
ind}_K^G(\sigma_i)$ and $\pi_i' = {^{\alpha}\!}(\pi_i)$, then we obtain
the unique supercuspidal $L$-packet $\{\pi_1,\pi_1',\pi_2,\pi_2'\}$ of
$G$ containing four elements.

For the ramified case, let $\widetilde{\sigma}$ be a very cuspidal
representation of $N_{\widetilde{G}}\widetilde{I}$ of level $l$ ($\geq
1$) and let $\widetilde{\pi}$ be the corresponding supercuspidal
representation of $\widetilde{G}$. Let $\sigma = {\rm
Res}_I(\widetilde{\sigma}).$ Then $\sigma = \sigma_1 \oplus \sigma_2$
for two irreducible representations $\sigma_i$ ($i =1,2$) of $I$ and
$\gamma$ conjugates one to the other, i.e., $\sigma_2 =
{^{\gamma}\!}\sigma_1$. Let $\pi_i = {\rm ind}_I^G(\sigma_i)$ and so
$\pi_2 = {^{\gamma}\!}\pi_1$. Then the restriction of $\widetilde{\pi}$
to $G$ breaks up into the direct sum of two irreducible supercuspidal
representations as ${\rm Res}_G(\widetilde{\pi}) = \pi_1 \oplus \pi_2.$
We call $\{\pi_1,\pi_2\}$ a {\it ramified supercuspidal $L$-packet} of
$G.$

To summarize, we have three kinds of supercuspidal $L$-packets for $G$
namely,
\begin{enumerate}
\renewcommand\labelenumi{(\roman{enumi})}
\leftskip .3pc
\item unramified supercuspidal $L$-packets $\{\pi,\pi'\}$;

\item the unique (unramified) supercuspidal $L$-packet
$\{\pi_1,\pi_1',\pi_2,\pi_2'\}$ of cardinality four;

\item ramified supercuspidal $L$-packets $\{\pi_1,\pi_2\}$.\vspace{-.5pc}
\end{enumerate}

\begin{propo} {\rm (Unramified supercuspidal $L$-packets of cardinality two)}
\label{prop:sc-unramified-2-sl2}$\left.\right.$\vspace{.5pc}

\noindent Consider an unramified supercuspidal $L$-packet $\{\pi,\pi'\}$
determined by a very cuspidal representation $\widetilde{\sigma}$ of
level $l$ of $Z\widetilde{K}$ as above. The conductors $c(\pi),c(\pi')$
are both equal to $2l$. The dimensions of the spaces $\pi_{\eta}^{K_m}$
and $(\pi')_{\eta}^{K_m}$ are as follows{\rm :}
\begin{enumerate}
\renewcommand\labelenumi{\rm (\roman{enumi})}
\leftskip .3pc
\item For any $\eta$ such that $\eta(-1) = \omega_{\pi}(-1)$ we have
\begin{equation*}
\hskip -1.25pc\pi^{K_{2l-1}}_{\eta} = \pi^{K_{2l-1}'}_{\eta} =
(\pi')^{K_{2l -1}}_{\eta} = (\pi')^{K_{2l-1}'}_{\eta} = (0).
\end{equation*}

\item Let $\eta(-1) = \omega_{\pi}(-1)$ and $c(\eta) \leq l.$ If $l$ is
odd then for all $m \geq 2l,$
%\begin{enumerate}
%\renewcommand\labelenumii{\rm (\alph{enumii})}
\leftskip .2pc
(a) $\dim (\pi^{K_m'}_{\eta}) = \dim
((\pi')^{K_m}_{\eta}) = 2 \left\lceil \frac{m-2l+1}{2}
\right\rceil,$\vspace{.3pc}

(b) $\dim (\pi^{K_m}_{\eta}) = \dim ((\pi')^{K_m'}_{\eta}) = 2
\left\lfloor \frac{m-2l+1}{2} \right\rfloor$.
%\end{enumerate}

\item Let $\eta(-1) = \omega_{\pi}(-1)$ and $c(\eta) \leq l.$ If $l$ is
even then for all $m \geq 2l,$
%\begin{enumerate}
%\renewcommand\labelenumii{\rm (\alph{enumii})}
\leftskip .2pc
(a) $\dim (\pi^{K_m}_{\eta}) =
\dim((\pi')^{K_m'}_{\eta}) = 2 \left\lceil \frac{m-2l+1}{2}
\right\rceil,$\vspace{.3pc}

(b) $\dim (\pi^{K_m'}_{\eta}) = \dim ((\pi')^{K_m}_{\eta}) = 2
\left\lfloor \frac{m-2l+1}{2} \right\rfloor.$
%\end{enumerate}
\end{enumerate}\vspace{-1.5pc}
\end{propo}

\begin{propo} {\rm (Test vectors for unramified
supercuspidal $L$-packets of cardinality two)}\label{prop:sc-unramified-2-sl2-test}
$\left.\right.$\vspace{.5pc}

\noindent Let $\widetilde{\sigma}$ be a very cuspidal representation of
$Z\widetilde{K}$ which determines an unramified supercuspidal $L$-packet
$\{\pi,\pi'\}$ as above. Assume that $\widetilde{\pi} = {\rm
ind}_{Z\widetilde {K}}^{\widetilde{G}}(\widetilde{\sigma})$ is realized
in its Kirillov model with respect to $\psi.$ Define two elements
$\phi_1$ and $\phi_{\epsilon}$ in the Kirillov model as follows{\rm :}
\begin{align*}
\phi_1(x) &= \left\{ \begin{array}{ll}
1, &\mbox{if $x \in (\mathcal{O}^{\times})^2,$}\\[.2pc]
0, &\mbox{if $x \notin (\mathcal{O}^{\times})^2,$}
\end{array} \right.\\[.2pc]
\phi_{\epsilon}(x) &= \widetilde{\pi}(\gamma)\phi_1.
\end{align*}
Let $\eta = \omega_{\widetilde{\pi}}$. We have
\begin{enumerate}
\renewcommand\labelenumi{\rm (\roman{enumi})}
\leftskip .3pc
\item ${\mathbb C}\phi_1 \oplus {\mathbb C}\phi_{\epsilon} =
\widetilde{\pi}^{K_{2l}}_{\eta}.$

\item If $l$ is even{\rm ,} then $\pi^{K_{2l}}_{\eta} =
\widetilde{\pi}^{K_{2l}}_{\eta}$. In addition{\rm ,} $\pi$ is
$\psi$-generic and any $\psi$-Whittaker\break\pagebreak \noindent functional is non-zero on
$\phi_1$ and vanishes on $\phi_{\epsilon}$. Furthermore{\rm ,} $\pi'$ is
not $\psi'$-generic for any character $\psi'$ of conductor
$\mathcal{O}$. It is however $\psi_{\varpi}$-generic{\rm ,} and any
$\psi_{\varpi}$-Whittaker functional is non-vanishing on
$\widetilde{\pi}(\alpha^{-1})\phi_1${\rm ,} which is a newform for
$\pi'.$

\item If $l$ is odd{\rm ,} then {\rm (ii)} holds with $\pi$ and $\pi'$
interchanged.
\end{enumerate}
\end{propo}

\begin{propo} {\rm (Unramified supercuspidal $L$-packet of cardinality four)}
\label{prop:sc-unramified-4-sl2}$\left.\right.$\vspace{.5pc}

\noindent Let $\widetilde{\sigma}$ denote a very cuspidal representation
of $Z\widetilde{K}$ of level $l=1$ such that ${\rm
Res}_K(\widetilde{\sigma}) = \sigma = \sigma_1 \oplus \sigma_2.$ Let
$\{\pi_1,\pi_1',\pi_2,\pi_2'\}$ be the corresponding $L$-packet of $G.$
Then $c(\pi_1) = c(\pi_1') = c(\pi_2) = c(\pi_2') = 2$. Moreover{\rm ,}
\begin{enumerate}
\renewcommand\labelenumi{\rm (\roman{enumi})}
\leftskip .2pc
\item Let $\eta$ be any character such that $\eta(-1) =
\omega_{\sigma}(-1).$ If $\pi$ denotes any representation in the
$L$-packet{\rm ,} then $\pi^{K_1}_{\eta} = \pi^{K_1'}_{\eta} = (0).$

\item Let $\eta$ be any character such that $\eta(-1) =
\omega_{\sigma}(-1)$ and $c(\eta) \leq 1$ then for all $m \geq 2$ we
have
%\begin{enumerate}
%\renewcommand\labelenumii{\rm (\alph{enumii})}
\leftskip .2pc
(a) $\dim ((\pi_1)^{K_m'}_{\eta}) = \dim
((\pi_2)^{K_m'}_{\eta}) = \dim ((\pi_1')^{K_m}_{\eta}) = \dim
((\pi_2')^{K_m}_{\eta}) = \left\lceil \frac{m-1}{2}
\right\rceil,$\vspace{.3pc}

(b) $\dim ((\pi_1)^{K_m}_{\eta}) = \dim ((\pi_1')^{K_m'}_{\eta}) =
\dim ((\pi_2)^{K_m}_{\eta}) = \dim ((\pi_2')^{K_m'}_{\eta}) =
\left \lfloor \frac{m-1}{2} \right \rfloor.$
%\end{enumerate}
\end{enumerate}
\end{propo}

\begin{propo} {\rm (Test vectors for unramified supercuspidal
$L$-packets of cardinality four)}\label{prop:sc-unramified-4-sl2-test}
$\left.\right.$\vspace{.5pc}

\noindent With notation as above let $\{\pi_1,\pi_1',\pi_2,\pi_2'\}$ be
the unramified supercuspidal $L$-packet of cardinality four. Let
$\overline{\psi}$ be the character of ${\mathbb F}_q$ obtained from
$\psi$ by identifying ${\mathbb F}_q$ with $\mathcal{P}^{-1}/
\mathcal{O}$. Without loss of generality assume that $\sigma_1$ is
$\overline{\psi}$-generic. Then
\begin{enumerate}
\renewcommand\labelenumi{\rm (\roman{enumi})}
\leftskip .2pc
\item $\pi_1'$ is $\psi$-generic{\rm ,} $\pi_1$ is
$\psi_{\varpi}$-generic{\rm ,} $\pi_2'$ is $\psi_{\epsilon}$-generic{\rm
,} and $\pi_2$ is $\psi_{\epsilon\varpi}$-generic.

\item Assume that $\widetilde{\pi}$ is realized in its $\psi$-Kirillov
model. The function $\phi_1$ of
Proposition~{\rm \ref{prop:sc-unramified-2-sl2-test}} is a newform for
$\pi_1'$. This further implies that $\widetilde{\pi}(\alpha)(\phi_1)$ is
a newform for $\pi_1,$ $\widetilde{\pi}(\gamma)(\phi_1)$ is a newform
for $\pi_2'$ and $\widetilde{\pi}(\alpha\gamma)(\phi_1)$ is a newform
for $\pi_2.$ Finally{\rm ,} each of these newforms is a test vector for an
appropriate Whittaker functional coming\break from {\rm (i)}.
\end{enumerate}
\end{propo}

\begin{propo} {\rm (Ramified supercuspidal $L$-packets)}
\label{prop:sc-ramified-sl2}$\left.\right.$\vspace{.5pc}

\noindent Let $\{\pi_1,\pi_2\}$ be a ramified supercuspidal $L$-packet
of level $l$ as above. Then $c(\pi_1) = c(\pi_2) = 2l+1$. Moreover{\rm ,}
\begin{enumerate}
\renewcommand\labelenumi{\rm (\roman{enumi})}
\leftskip .2pc
\item For any character $\eta$ of $F^*$ such that $\eta(-1) =
\omega_{\sigma}(-1)$ we have $(\pi_1)^{K_{2l}}_{\eta} =
(\pi_2)^{K_{2l}}_{\eta} = (\pi_1)^{K_{2l}'}_{\eta} =
(\pi_2)^{K_{2l}'}_{\eta} = (0).$

\item Let $\eta(-1) = \omega_\sigma(-1)$ and $c(\eta) \leq l.$ For all
$m \geq 2l+1$ we have $\dim ((\pi_1)^{K_m}_{\eta}) = \dim
((\pi_2)^{K_m}_{\eta}) = \dim ((\pi_1)^{K_m'}_{\eta}) = \dim
((\pi_2)^{K_m'}_{\eta}) = m-2l.$
\end{enumerate}
\end{propo}

\begin{propo} {\rm (Test vectors for ramified supercuspidal $L$-packets)}
\label{prop:sc-ramified-sl2-test}$\left.\right.$\vspace{.5pc}

\noindent Let $\{\pi_1,\pi_2\}$ be a ramified supercuspidal $L$-packet
coming from a very cuspidal representation $\widetilde{\sigma}$ of
$N_{\widetilde{G}}(\widetilde{I})$ of level $l\geq 1.$ One and only one
of the $\pi_i$ is $\psi$-generic{\rm ,} say $\pi_1$. Then $\pi_2$ is
$\psi_{\epsilon}$-generic. Let $\eta = \omega_{\sigma}.$ If $\phi_1$ and
$\phi_{\epsilon}$ have the same meaning as in
Proposition~{\rm \ref{prop:sc-unramified-2-sl2-test}} {\rm (}assuming that
$\widetilde{\pi}$ is realized in its Kirillov model{\rm ),} we have
\begin{enumerate}
\renewcommand\labelenumi{\rm (\roman{enumi})}
\leftskip .2pc
\item $(\pi_1)^{K_{2l+1}}_{\eta} = {\mathbb C}\phi_1$ and
$(\pi_2)^{K_{2l+1}}_{\eta} = {\mathbb C}\phi_{\epsilon}.$

\item Any $\psi$-Whittaker functional is non-zero on $\phi_1$ and
similarly any $\psi_{\epsilon}$-Whittaker functional is non-zero on
$\phi_{\epsilon}$.
\end{enumerate}
\end{propo}

\subsection{\it Comparison of conductor with other invariants}
\label{sec:comparison}

\setcounter{theore}{0}
\begin{theor}[\!]\label{thm:conductor}
Let $\pi$ be an irreducible admissible representation of $G = SL_2(F).$
Let $\widetilde{\pi}$ be a representation of $\widetilde{G} = GL_2(F)$
whose restriction to $G$ contains $\pi.$ Assume that $\widetilde{\pi}$
is minimal{\rm ,} i.e.{\rm ,} $c(\widetilde{\pi}\otimes \chi) \geq
c(\widetilde{\pi})$ for all characters $\chi$ of $F^*.$ Then
\begin{equation*}
c(\pi) = c(\widetilde{\pi}).
\end{equation*}
\end{theor}

The next theorem relates the conductor of a representation $\pi$ of $G$
with the depth (see~\cite{moy-prasad1}) $\rho(\pi)$ of $\pi$
(cf.~\cite{josh-raghuram1}).

\begin{theor}[(Relation between conductor and
depth)] \label{thm:conductor-depth-sl2}
Let $\pi$ be an irreducible representation of $G.$ Let $\rho(\pi)$ be
the depth of $\pi.$
\begin{enumerate}
\renewcommand\labelenumi{\rm (\roman{enumi})}
\leftskip .2pc
\item If $\pi$ is any subquotient of a principal series representation
$\pi(\chi)${\rm ,} then
\begin{equation*}
\hskip -1.25pc\rho(\pi) = {\rm max}\{c(\pi)-1,0\}.
\end{equation*}
\item If $\pi$ is an irreducible supercuspidal representation{\rm ,} then
\begin{equation*}
\hskip -1.25pc\rho(\pi) = {\rm max}\left\{\frac{c(\pi)-2}{2},0 \right\}.
\end{equation*}
\end{enumerate}
\end{theor}

\section{Newforms for ${\hbox{\bfseries\itshape U}\hbox{\bf \,(1,1)}}$}\label{sec:unitary}

\subsection{\it Definitions and preliminary remarks}\label{sec:defns-u11}

We now define the basic filtration subgroups of $\bar{G}$ as we did for
$G$ in~\S\ref{sec:newforms}. Let $\bar{K}_{-1}=\bar{G}$, $\bar{K}_0 =
\bar{K}$, the standard hyperspecial subgroup of $\bar{G}$, and for
$m\geq 1$,
\begin{equation*}
\bar{K}_m = \left\{ \left(\begin{array}{cc}
a &b\\
c &d\end{array}\right) \in \bar{K}:
c \equiv 0 \pmod{\mathcal{P}_E^m} \right\}.
\end{equation*}
We let $\bar{K}'_m = \alpha^{-1}\bar{K}_m\alpha$.

Let $(\bar{\pi},V)$ be an admissible representation of $\bar{G}$ such
that $\bar{Z}$ acts by scalars on $V$. Let $\bar{\eta}$ be a character
of $\mathcal{O}_E^{\times}$ such that $\bar{\eta}|_{E^1} =
\omega_{\bar{\pi}}$ (where we have identified $\bar{Z}$\break with $E^1$).

For any such character $\bar{\eta}$ and any subgroup $\bar{H}$ of
$\bar{G}$ we define
\begin{equation*}
\bar{\pi}_{\bar{\eta}}^{\bar{H}} := \left\{v \in V :
\bar{\pi}\left(\left(\begin{array}{cc}
a &b\\
c &d\end{array}\right)\right)v = \bar{\eta}(d)v, \quad
\forall \left(\begin{array}{cc}
a &b\\
c &d\end{array}\right) \in \bar{H} \right\}.
\end{equation*}
We define the {\it $\bar{\eta}$-conductor} $c_{\bar{\eta}}(\bar{\pi})$
of $\bar{\pi}$ to be
\begin{equation*}
c_{\bar{\eta}}(\bar{\pi}) = \min \{m :
\bar{\pi}_{\bar{\eta}}^{\bar{K}_m} \not= (0)\quad \mbox{or}\quad
\bar{\pi}_{\bar{\eta}}^{\bar{K}'_m} \not= (0)\}.
\end{equation*}
We define the {\it conductor} $c(\bar{\pi})$ of $\bar{\pi}$ as
\setcounter{equation}{0}
\begin{equation}
c(\bar{\pi}) = \min \{c_{\bar{\eta}}(\bar{\pi}) : \bar{\eta}|_{E^1} =
\omega_{\bar{\pi}}\}.
\end{equation}
If $\bar\eta$ is such that $c_{\bar\eta}(\bar\pi) = c(\bar\pi)$ and
$\bar\pi_{\bar\eta}^{\bar{K}_{c(\bar\pi)}}\neq (0)$
(resp.~$\bar\pi_{\bar\eta}^{\bar{K}'_{c(\bar\pi)}}\neq (0)$), then we
call $\bar\pi_{\bar\eta}^{\bar{K}_{c(\bar\pi)}}$
(resp.~$\bar\pi_{\bar\eta}^{\bar{K}'_{c(\bar\pi)}}$) a \textit{space of
newforms} of $\bar\pi$. In this case, we refer to a non-zero element of
$\bar\pi_{\bar\eta}^{\bar{K}_{c(\bar\pi)}}$ or
$\bar\pi_{\bar\eta}^{\bar{K}'_{c(\bar\pi)}}$ as a \textit{newform} of
$\bar\pi$.

In this section, we will compute the dimension of
$\bar{\pi}_{\bar{\eta}}^{\bar{K}_m}$ for every irreducible admissible
infinite-dimensional representation $\bar{\pi}$ of $\bar{G}$ and every
character $\bar{\eta}$ such that $c_{\bar{\eta}}(\bar{\pi}) =
c(\bar{\pi})$.

We will often make use of the following fact. Let $\pi$ be the
restriction of $\bar{\pi}$ to $G$, and let $\eta =
\bar{\eta}|_{\mathcal{O}_F^{\times}}.$ By definition, the group $K_m$
acts on $\bar{\pi}_\eta^{K_m}$ via the character $\eta$, hence via
$\bar{\eta}$. Also, $\bar{Z}$ acts on $\bar{\pi}_\eta^{K_m}$ via the
character $\omega_{\bar{\pi}}$, hence via $\bar{\eta}$ since
$\bar{\eta}|_{E^1} = \omega_{\bar{\pi}}$. Thus any
$\left(\begin{smallmatrix} a &b\\ c &d \end{smallmatrix}\right)
\in\bar{Z}K_m$ acts on $\bar{\pi}_\eta^{K_m}$ by multiplication by
$\bar{\eta}(d)$. In the light of (\ref{eq:filtration}),
\begin{align}\label{eq:index}
\bar{K}_m/\bar{Z}K_m &= \bar{T}_0/\bar{Z}T_0 \simeq
\mathcal{O}_E^{\times}/E^1\mathcal{O}_F^{\times} \simeq
N_{E/F}(\mathcal{O}_E^{\times})/N_{E/F}(E^1\mathcal{O}_F^{\times})\nonumber\\[.2pc]
&= \mathcal{O}_F^{\times}/(\mathcal{O}_F^{\times})^2.
\end{align}
We may therefore take $1$ and $\theta = \left(\begin{smallmatrix}
\epsilon_E &0\\ 0 &^{s}\epsilon_E^{-1}\end{smallmatrix}\right)$ as
coset representatives for $\bar{K}_m/\bar{Z}K_m$. Hence if
$v\in\bar{\pi}_\eta^{K_m}$, then
$v\in\bar{\pi}_{\bar{\eta}}^{\bar{K}_m}$ if and only if
$\bar{\pi}(\theta)v = \bar{\eta}( ^s\epsilon_E^{-1})v$, i.e.,
\begin{equation}\label{eq:theta}
\bar{\pi}_{\bar{\eta}}^{\bar{K}_m} = \left\{ v\in \bar{\pi}_{\eta}^{K_m}
: \bar{\pi}(\theta)v=\bar{\eta}(^s\epsilon_E^{-1})v\right\}.
\end{equation}

\vspace{.2pc}
\subsection{\it Principal series representations}\label{sec:principal-u11}\vspace{-.2pc}

Let $\bar{\chi}$ be a character of $E^*$. Let $\bar{\pi}(\bar{\chi})$
denote the principal series $\text{Ind}_{\bar{B}}^{\bar{G}}
(\bar{\chi})$. According to~\cite[\S11.1]{rogawski},
$\bar{\pi}(\bar{\chi})$ is irreducible except in the cases
\begin{enumerate}
\renewcommand\labelenumi{(\roman{enumi})}
\leftskip .2pc
\item $\bar{\chi}|_{F^*} = |\cdot |_F^{\pm}$,\vspace{.2pc}

\item $\bar{\chi}|_{F^*} = \omega_{E/F}$.
\end{enumerate}

In case (i), let $\mu$ be the character of $E^1$ defined by $\mu (a/^s
a) = \bar\chi |\cdot |_F^\mp$. Then $\bar{\pi}(\bar{\chi})$ has two
Jordan--H\"{o}lder constituents, namely the one-dimensional
representation $\xi = \mu\circ\det$ and a square integrable
representation $\text{St}(\xi)$. In case (ii), $\bar{\pi}(\bar{\chi})$
is the direct sum of two irreducible representations
$\bar{\pi}^1(\bar{\chi})$ and $\bar{\pi}^2(\bar{\chi})$, which together
form an $L$-packet of $\bar{G}$. We distinguish
$\bar{\pi}^1(\bar{\chi})$ from $\bar{\pi}^2(\bar{\chi})$ by defining
$\bar{\pi}^1(\bar{\chi})$ to be the summand that has a $K$-spherical
vector, hence $\bar{\pi}^i(\bar{\chi})|_G = \pi_E^i.$

Let $\chi = \bar{\chi}|_{F^*}$. Then the restriction of
$\bar{\pi}(\bar{\chi})$ to $G$ is isomorphic to $\pi (\chi)$. It is
easily seen then that the restriction to $G$ of any irreducible
constituent of $\bar{\pi}(\bar{\chi})$ is itself irreducible unless
$\chi$ is the character corresponding to some ramified quadratic
extension $E'/F$. In this case $\bar{\pi}(\bar{\chi})|_G$ decomposes as
the direct sum $\pi_{E'}^1 \oplus \pi_{E'}^2$. We now compute the
conductors of the representations in the principal series of $\bar{G}$.

\setcounter{theore}{0}
\begin{theor}[(Conductors for principal series
representations)] \label{thm:ps-u11}
Let $\bar{\chi}$ be a character of $E^*$. Suppose that $\bar{\pi}$ is an
irreducible constituent of the principal series $\bar{\pi}(\bar{\chi})$.
\begin{enumerate}
\renewcommand\labelenumi{\rm (\roman{enumi})}
\leftskip .2pc
\item If $\bar{\eta}$ is a character of $\mathcal{O}_E^\times$ with
$\bar{\eta}|_{E^1} = \omega_{\bar{\pi}}${\rm ,} then
$c_{\bar{\eta}}(\bar{\pi}) = c(\bar{\pi})$ if and only if
$\bar{\eta}=\bar{\chi}|_{\mathcal{O}_E^{\times}}$ or
$^s\!\bar{\chi}^{-1}|_{\mathcal{O}_E^{\times}}$. Moreover{\rm ,}
\begin{align*}
\hskip -1.25pc c(\bar{\pi}) = \left\{\begin{array}{ll}
c(\bar{\chi}|_{F^*}), &\mbox{if $\bar{\pi}\neq {\rm St}(\xi)$},\\[.2pc]
1, &\mbox{if $\bar{\pi} = {\rm St}(\xi)$}.
\end{array}\right.
\end{align*}

\item Suppose $\bar{\eta}$ is as above.
%\begin{enumerate}
%\renewcommand\labelenumii{\rm (\alph{enumii})}
\leftskip .2pc
(a) If $\bar{\pi} = \bar{\pi}(\bar{\chi})${\rm ,}
$\bar{\chi}$ is ramified{\rm ,} and
$\bar{\chi}|_{\mathcal{O}_E^{\times}} =
^{s}\bar{\chi}^{-1}|_{\mathcal{O}_E^{\times}}${\rm ,} then
\begin{equation*}
\hskip -2.5pc \dim (\bar{\pi}_{\bar{\eta}}^{\bar{K}_m}) = \left\{ \begin{array}{ll}
0, & \mbox{if $m=0$},\\[.2pc]
m+1, & \mbox{if $m>0$}. \end{array} \right.
\end{equation*}

(b) For $\bar{\pi}=\bar{\pi}^1(\bar{\chi}),
\bar{\pi}^2(\bar{\chi})${\rm ,} we have
\begin{align*}
\hskip -2.5pc \dim\left(\bar{\pi}^1(\bar{\chi})_{\bar{\eta}}^{\bar{K}_m}\right) =
\dim\left(\bar{\pi}^2(\bar{\chi})_{\bar{\eta}}^{\bar{K}'_m}\right) &=
\left\lceil\frac{m+1}{2}\right\rceil,\\[.2pc]
\hskip -2.5pc \dim\left(\bar{\pi}^2(\bar{\chi})_{\bar{\eta}}^{\bar{K}_m}\right) =
\dim\left(\bar{\pi}^1(\bar{\chi})_{\bar{\eta}}^{\bar{K}'_m}\right) &=
\left\lceil\frac{m}{2}\right\rceil.
\end{align*}

(c) In all other cases{\rm ,}
\begin{equation*}
\hskip -2.5pc \dim\bar{\pi}_{\bar{\eta}}^{\bar{K}_m} = \max\{ m-c(\bar{\pi})+1,0\}.
\end{equation*}
%\end{enumerate}
\end{enumerate}
\end{theor}

\begin{proof}
We may assume without loss of generality that $\bar{\chi}$ is chosen so
that $\bar{\pi}$ is a subrepresentation of $\bar{\pi}(\bar{\chi})$. Let
$\pi$ be the restriction of $\bar{\pi}$ to $G$. Let $\bar{\eta}$ be any
character of $\mathcal{O}_E^{\times}$ with $\bar{\eta}|_{E^1} =
\omega_{\bar{\pi}}$. Let $\eta=\bar{\eta}|_{\mathcal{O}_F^{\times}}$.
Since $\bar{\pi}_{\bar{\eta}}^{\bar{K}_m}\subset \pi_{\eta}^{K_m}$,
\begin{equation*}
c_{\bar{\eta}}(\bar{\pi})\geq c_{\eta}(\pi)\geq c(\pi).
\end{equation*}
We claim that $c_{\bar{\eta}}(\bar{\pi}) = c(\pi)$ precisely for
$\bar{\eta} = \bar{\chi}$ or $^s\bar{\chi}^{-1}$.  The first part of (i)
follows immediately from this claim, and the second follows from this
together with the conductor calculations in~\S\ref{sec:principal-sl2}.

Let $c=c(\pi)$. The only $\eta$ such that $c_\eta (\pi) = c (\pi)$ are
$\bar{\chi}|_{\mathcal{O}_F^{\times}}$ and
$\bar{\chi}^{-1}|_{\mathcal{O}_F^{\times}}$. Hence we cannot have
$c_{\bar{\eta}}(\bar{\pi}) = c(\pi)$ unless $\bar{\eta}$ equals
$\bar{\chi}^\pm$ on $\mathcal{O}_F^{\times}$. We first prove that
$\bar{\pi}_{\bar{\eta}}^{\bar{K}_c}\neq (0)$ if and only if $\bar{\eta}
= \bar{\chi}|_{\mathcal{O}_E^{\times}}$ or
$^s\!\bar{\chi}^{-1}|_{\mathcal{O}_E^{\times}}$ in the case where
$\bar{\eta}|_{\mathcal{O}_F^{\times}} = \eta =
\bar{\chi}|_{\mathcal{O}_F^{\times}}$ and $\bar{\pi}\neq\bar{\pi}^2(\bar{\chi})$.

Since $\pi_{\eta}^{K_c}$ is contained in the restriction of
$\bar{\pi}(\bar{\chi})$ to $G$, which is isomorphic to $\pi
(\bar{\chi}|_{F^*})$, it is an easy consequence of the proofs of the
statements in~\S\ref{sec:principal-sl2} (see \cite{josh-raghuram2}) that
\begin{align*}
\pi_{\eta}^{K_c} = \left\{ \begin{array}{ll}
\mathbb{C}\bar{f}_w, &\mbox{if $\bar{\chi}^2|_{\mathcal{O}_F^{\times}}
\neq {\bf 1}$,}\\[.2pc]
\mathbb{C}\bar{f}_w + \mathbb{C}\bar{f}_1, &\mbox{if
$\bar{\chi}^2|_{\mathcal{O}_F^{\times}} = {\bf 1}$,}
\end{array} \right.
\end{align*}
where
\begin{align*}
\bar{f}_w (g) &= \left\{ \begin{array}{ll}
0, &\mbox{if $g\notin\bar{B}wK_c$,}\\[.2pc]
\bar{\chi}(t)|t|_E^{1/2}\eta (d) = \bar{\chi}(td)|t|_E^{1/2}, &\mbox{if
$g= \left(\begin{array}{cc}
t &*\\
0 &^{s}t^{-1} \end{array}\right) w \left(\begin{array}{cc}
a &b\\
c &d\end{array}\right)$,} \end{array} \right.\\[.2pc]
\bar{f}_1 (g) &= \left\{ \begin{array}{ll}
0, &\mbox{if $g\notin\bar{B}K_c$,}\\[.2pc]
\bar{\chi}(t)|t|_E^{1/2}\eta (d) = \bar{\chi}(td)|t|_E^{1/2}, &\mbox{if
$g = \left(\begin{array}{cc}
t &*\\
0 &^{s}t^{-1} \end{array}\right) \left(\begin{array}{cc}
a &b\\
c &d\end{array}\right)$.} \end{array} \right.
\end{align*}

We now determine when $\bar{f}_w,\bar{f}_1$ lie in
$\bar{\pi}_{\bar{\eta}}^{\bar{K}_c}$. In the light of (\ref{eq:theta}),
this reduces to verifying whether $\bar{\pi}(\theta)$ acts as the scalar
$\bar{\eta}( ^s\epsilon_E^{-1})$ on these vectors. It is easily checked
that
\begin{align*}
\bar{\pi}(\theta)\bar{f}_w &= \bar{\chi}(^s\epsilon_E^{-1})\bar{f}_w,\\[.2pc]
\bar{\pi}(\theta)\bar{f}_1 &= \bar{\chi}(\epsilon_E)\bar{f}_1.
\end{align*}

Hence $\bar{f}_w\in\bar{\pi}_{\bar{\eta}}^{\bar{K}_c}$ if and only if
$\bar{\eta}(^s\epsilon_E) = \bar{\chi}(^s\epsilon_E)$. This is
equivalent to $\bar{\eta}= \bar{\chi}|_{\mathcal{O}_E^{\times}}$ since
$\bar{\eta}$ and $\bar{\chi}$ already agree on $\mathcal{O}_F^{\times}$
and $E^1$ (by assumption) and since $^s\epsilon_E$ is a representative
for the non-trivial coset in
$\mathcal{O}_E^\times/E^1\mathcal{O}_F^\times$. Similarly, if
$\bar{\chi}^2|_{\mathcal{O}_F^{\times}} = {\bf 1}$, then
$\bar{f}_1\in\bar{\pi}_{\bar{\eta}}^{\bar{K}_c}$ if and only if
$\bar{\eta}(\epsilon_E) =\, ^s\bar{\chi}^{-1}(\epsilon_E)$, which is
equivalent to $\bar{\eta}=\,
^s\bar{\chi}^{-1}|_{\mathcal{O}_E^{\times}}$ since $\bar{\eta}$ and
$^s\bar{\chi}^{-1}$ already agree on $\mathcal{O}_F^{\times}$ and $E^1$
and since the non-trivial coset in
$\mathcal{O}_E^\times/E^1\mathcal{O}_F^\times$ is represented by
$\epsilon_E$. Summarizing, we have that when
$\bar{\pi}\neq\bar{\pi}^2(\bar{\chi})$ and $\eta =\bar{\chi}
|_{\mathcal{O}_F^{\times}}$, $\bar{\pi}_{\bar{\eta}}^{\bar{K}_c}\neq
(0)$ if and only if $\bar{\eta} = \bar{\chi}|_{\mathcal{O}_E^{\times}}$
or $^s\!\bar{\chi}^{-1}|_{\mathcal{O}_E^{\times}}$.

On the other hand, if $\eta =
\bar{\chi}^{-1}|_{\mathcal{O}_F^{\times}}$, note that we may exchange
$\bar{\chi}$ and $^s\bar{\chi}^{-1}$ in the above proof since
$\bar{\pi}(\bar{\chi})$ and $\bar{\pi}(^s\bar{\chi}^{-1})$ have the same
constituents. (Of course, exchanging $\bar{\chi}$ and
$^s\bar{\chi}^{-1}$ may make our assumption that $\bar{\pi}$ is a
subrepresentation of $\bar{\pi}(\bar{\chi})$ false. The only case in
which this matters, however, is when $\bar{\pi} = \text{St}(\xi)$, and
in this case we are already done since $\bar{\chi} = \bar{\chi}^{-1}$ on
$\mathcal{O}_F^{\times}$.) Then carrying out the proof \textit{mutatis
mutandis}, we obtain again that $\bar{\pi}_{\bar{\eta}}^{\bar{K}_c}\neq
(0)$ if and only if $\bar{\eta} = \bar{\chi}$ or $^s\bar{\chi}^{-1}$.
This establishes our claim if $\bar{\pi}$ is in a singleton $L$-packet
since for all $m\geq 0$,
\begin{equation*}
\dim (\bar{\pi}_{\bar{\eta}}^{\bar{K}'_m}) =
\dim (\, ^\alpha \bar{\pi}_{\bar{\eta}}^{\bar{K}_m})  =
\dim (\bar{\pi}_{\bar{\eta}}^{\bar{K}_m}).
\end{equation*}

Finally, suppose that $\bar{\pi} = \bar{\pi}^1(\bar{\chi})$. By the
above, $\left(\bar{\pi}^1(\bar{\eta})\right) _{\bar{\eta}}^{\bar{K}_0}
\neq (0)$ if and only if $\bar{\eta} = \bar{\chi}$ or
$^s\bar{\chi}^{-1}$. Also, if $\bar{\eta}$ is any character of
$\mathcal{O}_E^{\times}$, then since
$\left(\bar{\pi}^1(\bar{\chi})\right) _{\bar{\eta}}^{\bar{K}'_0} =
\left( ^\alpha\bar{\pi}^1(\bar{\chi})\right) _{\bar{\eta}}^{\bar{K}_0}$
and $^\alpha\bar{\pi}^1(\bar{\chi})\cong\bar{\pi}^2(\bar{\chi})$, we
have that
\begin{equation*}
\dim (\bar{\pi}^1(\bar{\chi})) _{\bar{\eta}}^{\bar{K}'_0}
= \dim (^\alpha\bar{\pi}^1(\bar{\chi}))
_{\bar{\eta}}^{\bar{K}_0} = \dim (\bar{\pi}^2(\bar{\chi}))
_{\bar{\eta}}^{\bar{K}_0}.
\end{equation*}
But $\left( \bar{\pi}^2(\bar{\chi})\right)
_{\bar{\eta}}^{\bar{K}_0}\subset\left( \bar{\pi}^2(\bar{\chi})\right)
_{\eta}^{K_0} = (0)$ by Proposition~\ref{prop:unramfified-sl2} since
$\text{Res}_G\,\bar{\pi}^2(\bar{\chi})\cong\pi_E^2$. Thus
$\dim\left(\bar{\pi}^1(\bar{\chi})\right) _{\bar{\eta}}^{\bar{K}'_0}=0$
so again $c_{\bar{\eta}}(\bar{\pi}) = 0 = c(\pi)$ precisely for
$\bar{\eta} = \bar{\chi}$ or $^s\bar{\chi}^{-1}$. Finally, conjugating
by $\alpha$ as above, one easily obtains the claim in the case
$\bar{\pi} = \bar{\pi}^2(\bar{\chi})$.

We now compute the dimensions of $\bar{\pi}_{\bar{\eta}}^{\bar{K}_m}$ to
prove (ii). Since $\bar{\pi} (\bar{\chi})$ and $\bar{\pi}
(^s\bar{\chi}^{-1})$ have the same irreducible constituents, we may
assume that $\bar{\eta} = \bar{\chi}|_{\mathcal{O}_E^{\times}}$. (As
above, the representations $\text{St}(\xi)$ present no problem here
since in this case $\bar{\chi} = \, ^s\bar{\chi}^{-1}$ on
$\mathcal{O}_E^{\times}$.)

If $\bar{\pi}\neq\bar{\pi}^2(\bar{\chi})$ the proof of (i) shows that
$\bar{\pi}_{\bar{\eta}}^{\bar{K}_c} = \pi_\eta^{K_c}$. Thus
$\dim\bar{\pi}_{\bar{\eta}}^{\bar{K}_c}$ is 1 if
$\bar{\chi}|_{\mathcal{O}_E^{\times}}\neq\,
^s\bar{\chi}^{-1}|_{\mathcal{O}_E^{\times}}$ and 2 if
$\bar{\chi}|_{\mathcal{O}_E^{\times}} = \,
^s\bar{\chi}^{-1}|_{\mathcal{O}_E^{\times}}$. The proof also shows that
\begin{align*}
\dim (\bar{\pi}^1(\bar{\chi}))_{\bar{\eta}}^{\bar{K}_0} =
\dim (\bar{\pi}^2(\bar{\chi}))_{\bar{\eta}}^{\bar{K}'_0} &=
1,\\[.2pc]
\dim (\bar{\pi}^2(\bar{\chi}))_{\bar{\eta}}^{\bar{K}_0} =
\dim (\bar{\pi}^1(\bar{\chi}))_{\bar{\eta}}^{\bar{K}'_0} &=
0.
\end{align*}
This shows that the formulae for the dimensions are valid when $m=c$.\pagebreak

Suppose that $m>c$. As with Theorem~5.3 of \cite{dipendra-raghuram}, it
follows from Lemma~3.2.1 and the following proofs in \S 3.2 of
\cite{josh-raghuram2} that $\pi_{\eta}^{K_m}$ is the direct sum of
$\pi_{\eta}^{K_c}$ together with certain two-dimensional spaces
$\pi_{\eta ,i}^{K_m}$ of the form
$\mathbb{C}\bar{f}_{i,1}+\mathbb{C}\bar{f}_{i,\epsilon}$ ($1\leq i\leq
m-c$), where
\begin{align*}
\hskip -4pc \bar{f}_{i,1} (g) &= \left\{ \begin{array}{l@{\quad}l}
0, &\mbox{if $g\notin\bar{B} \left(\begin{array}{cc}
1 &0\\
\varpi^m &1 \end{array}\right) K_m$,}\\[.3pc]
\bar{\chi}(t)|t|_E^{1/2}\eta (d) = \bar{\chi}(td)|t|_E^{1/2},
&\mbox{if $g = \left(\begin{array}{cc}
t &*\\
0 &^st^{-1} \end{array}\right) \left(\begin{array}{cc}
1 &0\\
\varpi^i &1 \end{array}\right) \left(\begin{array}{cc}
a &b\\
c &d\end{array}\right)$,} \end{array} \right.\\[.3pc]
\hskip -4pc \bar{f}_{i,\epsilon} (g) &= \left\{ \begin{array}{l@{\quad}l}
0, &\mbox{if $g\notin\bar{B} \left(\begin{array}{cc}
1 &0\\
\varpi^m\epsilon_F &1 \end{array}\right) K_m$,}\\[.3pc]
\bar{\chi}(t)|t|_E^{1/2}\eta (d) = \bar{\chi}(td)|t|_E^{1/2},
&\mbox{if $g= \left(\begin{array}{cc}
t &*\\
0 &^st^{-1} \end{array}\right) \left(\begin{array}{cc}
1 &0\\
\varpi^i\epsilon_F &1 \end{array}\right) \left(\begin{array}{cc}
a &b\\
c &d\end{array}\right)$.} \end{array} \right.
\end{align*}
We will now verify that whenever $\pi_{\eta ,i}^{K_m} \subset
\pi_\eta^{K_m}$,
\begin{enumerate}
\leftskip .2pc
\item $\pi_{\eta ,i}^{K_m}$ is $\bar{K}_m$-stable, and

\item the subspace of $\pi_{\eta ,i}^{K_m}$ on which $\bar{K}_m$ acts
via the character $\bar{\eta}$ is one-dimensional.
\end{enumerate}
If this holds, then
\begin{equation*}
\dim \bar{\pi}_{\bar{\eta}}^{\bar{K}_m}- \dim
\bar{\pi}_{\bar{\eta}}^{\bar{K}_c} = \frac{1}{2} (\dim
\pi_{\eta}^{K_m}- \dim\pi_{\eta}^{K_c}),
\end{equation*}
and the formulae for the dimension of
$\bar{\pi}_{\bar{\eta}}^{\bar{K}_m}$ follow easily from this equation
and the dimension results of~\S\ref{sec:principal-sl2}. The dimension of
$\bar{\pi}_{\bar{\eta}}^{\bar{K}'_m}$ is computed analogously.

We now show (1) and (2). By (\ref{eq:theta}), this reduces to showing
that $\pi_{\eta ,i}^{K_m}$ is $\theta$-stable, and that the subspace of
$\pi_{\eta ,i}^{K_m}$ on which $\theta$ acts as the scalar $ \bar{\eta}(
^s\epsilon_E^{-1})$ is one-dimensional. Let $\delta$ be either 1 or
$\epsilon$. Then
\begin{equation*}
(\bar{\pi}(\theta)\bar{f}_{i,\delta})(g) = \bar{f}_{i,\delta}(g\theta)
= \bar{f}_{i,\delta} (\theta\ (\theta^{-1}g\theta))
= \bar{\chi}(\epsilon_E)\bar{f}_{i,\delta}(\theta^{-1} g\theta).
\end{equation*}
If $\delta=1$, this is non-zero if and only if $\theta ^{-1} g\theta\in
\bar{B} \left(\begin{smallmatrix} 1 &0\\ \varpi^i &1
\end{smallmatrix}\right) K_m$, i.e., if and only if $g\in \bar{B}
\left(\begin{smallmatrix} 1 &0\\ \varpi^i\epsilon_F &1
\end{smallmatrix}\right) K_m$. This together with the fact that
$\bar{\pi}(\theta)\bar{f}_{i,1}\in\pi_\eta^{K_m}$ implies that
$\bar{\pi}(\theta)\bar{f}_{i,1}$ is a multiple of
$\bar{f}_{i,\epsilon}$. The exact multiple is determined by evaluating
\begin{align*}
\left(\bar{\pi}(\theta)\bar{f}_{i,1}\right) \left( \left(\begin{array}{cc}
1 &0\\
\varpi^i\epsilon_F &1 \end{array}\right) \right)
&= \bar{f}_{i,1}\left( \left(\begin{array}{cc}
1 &0\\
\varpi^i\epsilon_F &1 \end{array}\right) \theta\right)\\
&= \bar{f}_{i,1}\left( \left(\begin{array}{cc}
^s\epsilon_E^{-1} &0\\
0 &\epsilon_E \end{array}\right) \left(\begin{array}{cc}
1 &0\\
\varpi^i &1 \end{array}\right) \left(\begin{array}{cc@{\ }}
\epsilon_F &0\\
0 &\epsilon_F^{-1} \end{array}\right) \right)\\
&= \bar{\chi}(^s\epsilon_E^{-1}\epsilon_F^{-1}).
\end{align*}
Thus
\begin{equation*}
\bar{\pi}(\theta)\bar{f}_{i,1} = \bar{\chi} (^s\epsilon_E^{-1}
\epsilon_F^{-1})\bar{f}_{i,\epsilon}.
\end{equation*}
Similarly,
\begin{equation*}
\bar{\pi}(\theta) \bar{f}_{i,\epsilon} = \bar{\chi} (\epsilon_E)
\bar{f}_{i,1}.
\end{equation*}

As claimed, $\theta$ stabilizes $\bar{\pi}_{\eta,i}^{K_m}$. Moreover,
the characteristic polynomial of $\theta$ on this two-dimensional space
is
\begin{equation*}
X^2 - \bar{\chi}(^s\epsilon_E^{-1}\epsilon_F^{-1})\bar{\chi}(\epsilon_E)
= X^2-\bar{\chi}( ^s\epsilon_E^{-1})^2.
\end{equation*}
The eigenvalues of $\theta$ on $\bar{\pi}_{\eta ,i}^{K_m}$ are therefore
$\pm \bar{\chi}(^s\epsilon_E^{-1}) = \pm \bar{\eta}(^s\epsilon_E^{-1})$.
It follows that the subspace of $\pi_{\eta ,i}^{K_m}$ on which $\theta$
acts as the scalar $\bar{\eta}( ^s\epsilon_E^{-1})$ is
one-dimensional.\hfill $\Box$
\end{proof}

Now suppose that $\bar{\pi}$ is an irreducible representation of
conductor $c$ in the principal series of $\bar{G}$ and that $\bar{\eta}$
is such that $c_{\bar{\eta}}(\bar{\pi}) = c$. We consider the effect of
the Whittaker functional $\Lambda_\psi$ given by
(\ref{eqn:whittaker-ps}) on $\bar{\pi}_{\bar{\eta}}^{\bar{K}_c},
\bar{\pi}_{\bar{\eta}}^{\bar{K}'_c}$.

\begin{propo} {\rm (Test vectors for principal series representations)}
\label{prop:ps-test-u11}$\left.\right.$\vspace{.5pc}

\noindent Suppose that $\bar{\pi}$ is an irreducible representation in
the principal series of $\bar{G}$. Let $\bar{\eta}$ be a character of
$\mathcal{O}_E^{\times}$ with $\bar{\eta}|_{E^1} = \omega_{\bar{\pi}}$
such that $c_{\bar{\eta}}(\bar{\pi}) = c (\bar{\pi})$. Let $\psi =
\psi_F$.
\begin{enumerate}
\renewcommand\labelenumi{\rm (\roman{enumi})}
\leftskip .5pc
\item If $\bar{\pi} = \bar{\pi}(\bar{\chi})${\rm ,} $\bar{\chi}$ is
ramified{\rm ,} and $\bar{\chi} = \,
^s\bar{\chi}^{-1}|_{\mathcal{O}_E^{\times}}${\rm ,} then $\bar{\pi}$ is
$\psi$-generic. Moreover{\rm ,} the space of vectors
$\bar{\pi}_{\bar{\eta}}^{\bar{K}_1}$ on which $\Lambda_\psi$ vanishes is
one-dimensional.

\item If $\bar{\pi} = \bar{\pi}^1(\bar{\chi})${\rm ,} then $\bar{\pi}$
is $\psi$-generic and $\Lambda_\psi$ is non-zero on the one-dimensional
space of newforms $\bar{\pi}_{\bar{\eta}}^{\bar{K}_0}$.

\item If $\bar{\pi} = \bar{\pi}^2(\bar{\chi})$, then $\bar{\pi}$ is not
$\psi$-generic{\rm ,} but it is $\psi_\varpi$-generic. Moreover{\rm ,}
$\Lambda_{\psi_\varpi}$ is non-zero on the one-dimensional space of
newforms $\bar{\pi}_{\bar{\eta}}^{\bar{K}'_0}$.

\item In all other cases{\rm ,} $\bar{\pi}$ is $\psi$-generic. In
addition{\rm ,} if $c=c(\bar{\pi})${\rm ,} then $\Lambda_\psi (v)\neq 0$
for any newform $v$ in $\bar{\pi}_{\bar{\eta}}^{\bar{K}_c}$.
\end{enumerate}\vspace{-1.5pc}
\end{propo}

\begin{proof}
Let $\pi$ be the restriction of $\bar{\pi}$ to $G$. Note that since $G$
and $\bar{G}$ have Borel subgroups with the same unipotent radical
(namely, $N$), the restriction of $\Lambda_\psi$ to any $\psi$-generic
component of $\pi$ is a non-zero $\psi$-Whittaker functional on that
component, while its restriction to any non-$\psi$-generic component is
$0$.

Let $c=c(\bar\pi)$. Assume we are in case (ii), (iii), or (iv). Let
$\bar{L}$ be either $\bar{K}_c$ or $\bar{K}'_c$, according to the case,
and let $L = \bar{L}\cap G$, i.e., $L$ is either $K_c$ or $K'_c$. Assume
that $v$ is a non-zero vector in $\bar{\pi}_{\bar{\eta}}^{\bar{L}}$. By
Theorem~\ref{thm:ps-u11}, the restriction of $\bar{\pi}$ to $G$ is
irreducible of conductor $c$, and $\bar{\pi}_{\bar{\eta}}^{\bar{L}} =
\pi_\eta^{L}$ is one-dimensional. The statements in each of these cases
now follow easily from the analogous results about $\pi$
in~\S\ref{sec:principal-sl2}.

Suppose now that $\bar{\pi} = \bar{\pi}(\bar{\chi})$ with $\bar{\chi}$
ramified and $\bar{\chi} =\,^s \bar{\chi}^{-1}$. Then $\pi$ has
conductor $c=1$ and $\bar{\pi}_{\bar{\eta}}^{\bar{K}_1} =
\pi_\eta^{K_1}$ has dimension 2.

If $\pi$ is irreducible, then $\pi$ is $\psi$-generic according to
Corollary~\ref{cor:test-ram-ps}. Also, according to the proof of
Theorem~\ref{thm:ps-u11} (and using its notation),
$\bar{\pi}_{\bar{\eta}}^{K_1} = \mathbb{C}f_1\oplus\mathbb{C}f_w$. It
follows from Corollary~\ref{cor:test-ram-ps} that $\Lambda_\psi
(f_1)\neq 0$. Since the image of $\Lambda_\psi$ has dimension 1,
$\Lambda_\psi$ must vanish on a one-dimensional subspace of
$\bar{\pi}_{\bar{\eta}}^{\bar{K}_1}$.

If $\pi$ is reducible, then as discussed in~\S\ref{sec:principal-sl2},
$\pi$ decomposes as the direct sum of two representations $\pi_1$ and
$\pi_2$. Moreover, only one of these representations, say $\pi_1$, is
$\psi$-generic by Corollary~\ref{cor:test-ramified-l-packets-sl2}. Then
$\Lambda_\psi$ vanishes on \hbox{$\left(\pi_2\right)_\eta^{K_1}\subset\pi_\eta^{K_1} =
\bar{\pi}_{\bar{\eta}}^{\bar{K}_1}$.} Moreover, by
Corollary~\ref{cor:test-ramified-l-packets-sl2}, \hbox{$\Lambda_\psi (v)\neq 0$} for all non-zero
\hbox{$v\in\left(\pi_1\right)_\eta^{K_1}\subset\pi_\eta^{K_1} =
\bar{\pi}_{\bar{\eta}}^{\bar{K}_1}$.} Hence, as in the preceding
paragraph, the subspace of $\bar{\pi}_{\bar{\eta}}^{\bar{K}_1}$ on which
$\Lambda_\psi$ vanishes is one-dimensional.\hfill $\Box$
\end{proof}

\vspace{.05pc}
\subsection{\it Supercuspidal representations}\label{sec:sc-u11}

We now consider the supercuspidal representations of $\bar{G}$. Let
$\bar{\pi}$ be such a representation. It is easily deduced from
analogous results on $\widetilde{G}$ and $G$ that $\bar{\pi}$ is
compactly induced from an irreducible representation of $\bar{K}$,
$\bar{K}'$, or $\bar{I}$. We will call $\bar{\pi}$ an \textit{unramified
{\rm (}resp. ramified{\rm )} supercuspidal representation} of $\bar{G}$ if its
restriction to $G$ contains an unramified (resp. ramified) supercuspidal
representation of $G.$

\subsubsection*{\it \!Ramified case.}

Suppose first that $\bar{\pi}$ is ramified. Let $\pi$ be the restriction
of $\bar{\pi}$ to $G$. Let $\pi_1$ be an irreducible component of the
restriction of $\bar{\pi}$ to $G$. Then $\pi_1$ is a ramified
supercuspidal representation of $G$. We extend $\pi_1$ to a
representation of $\bar{Z}G$ via the central character
$\omega_{\bar{\pi}}$, also denoted by $\pi_1$. Then $\bar{\pi}$ is
contained in $\text{ind}_{\bar{Z}G}^{\bar{G}}\pi_1$, and the restriction
of $\text{ind}_{\bar{Z}G}^{\bar{G}}\pi_1$ to $\bar{Z}G$ is
$\pi_1\oplus\,^\theta\pi_1$. But conjugation by $\theta$ and $\gamma$
have the same effect on $G$ so, by the discussion in the beginning
of~\S\ref{sec:sc-sl2}, $\pi_1$ and $\pi_2=\,^\theta\pi_1$ comprise an
$L$-packet for $G$. Since $\pi_1 \not\cong\,^\theta\pi_1$,
$\text{ind}_{\bar{Z}G}^{\bar{G}}\pi_1$ is irreducible and hence equal to
$\bar{\pi}$. Thus $\text{Res}_G\,\bar{\pi} = \pi_1\oplus\pi_2$, where\break
$\pi_2\cong\,^\theta\pi_1$.

From Theorem~\ref{thm:conductor-depth-sl2}, we see that the conductor of
both $\pi_1$ and $\pi_2$ is $2\rho + 2$, where $\rho$ is the depth of
both $\pi_1$ and $\pi_2$. We note that the depth of a twist of
$\bar{\pi}$ is no less than $\rho$. To see this, let $x$ be a point in
the Bruhat--Tits building of $\bar{G}$ (which is the same as that of $G$)
and let $r$ be a non-negative real number. Then any vector in the twist
of $\bar{\pi}$ that is fixed by $\bar{G}_{x,r+}$ is fixed by $G_{x,r+}$
since $G_{x,r+}\subset\bar{G}_{x,r+}$ (see~\cite{moy-prasad1}). It
follows that the depth of the twist of $\bar{\pi}$ is no less than the
depth of its restriction to $G$. But this restriction is $\pi$, which
has depth equal to $\rho$.

On the other hand, we may select a character $\chi$ of $\bar{G}$ such
that $\chi^{-2} = \omega_{\bar{\pi}}$ on $E^1\cap (1+\mathcal{P}_E)$
(viewed as a subgroup of $\bar{Z}$). If $\bar{\pi}' = \bar{\pi}\otimes
\chi$, then $\omega_{\bar{\pi}\otimes\chi}$ is trivial on $E^1\cap
(1+\mathcal{P}_E)$, and it is easily seen that $\rho (\bar{\pi}') =
\rho$. Define $\rho_0(\bar{\pi}) = \min\{ \rho (\bar{\pi}\otimes
\chi)\}$ as $\chi$ ranges over all characters of $\bar{G}$. Then we have
$\rho_0(\bar{\pi})=\rho$.

\setcounter{theore}{0}
\begin{theor}[(Ramified supercuspidal
representation)] \label{thm:sc-ram-u11}
Let $(\bar{\pi},V)$ be a ramified supercuspidal representation of
$\bar{G}$. Let $\bar{\eta}$ be any character of $\mathcal{O}_E^{\times}$
with $\bar{\eta}|_{E^1} = \omega_{\bar{\pi}}$ and
$c(\bar{\eta}|_{\mathcal{O}_F^{\times}}) \leq \rho_0(\bar{\pi})+ 1/2$.
Then we have $c(\bar{\pi}) = c_{\bar{\eta}}(\bar{\pi}) =
2\rho_0(\bar{\pi}) + 2$ and
\begin{equation*}
\dim (\bar{\pi}_{\bar{\eta}}^{\bar{K}_m}) = \max\left\{ m-c(\bar{\pi}) +
1, 0\right\}.
\end{equation*}
\end{theor}

\vspace{.1pc}
\begin{proof}
Let $\pi$ be the restriction of $\bar{\pi}$ to $G$. Set
$c=2\rho_0(\bar{\pi}) +2$ and
$\eta=\bar{\eta}|_{\mathcal{O}_F^{\times}}$. As discussed above, the
restriction of $\bar{\pi}$ to $G$ is the direct sum of two ramified
supercuspidal representations $\pi_1,\pi_2$, each of conductor $c$. By
Proposition~\ref{prop:sc-ramified-sl2}, $\dim (\pi_1)_\eta^{K_m} = \dim
(\pi_2)_\eta^{K_m}$ is non-zero if and only if $m\geq c$. Hence if $m<
c$, $\dim\bar{\pi}_{\bar{\eta}}^{\bar{K}_m}=0$ since
\begin{equation*}\bar{\pi}_{\bar{\eta}}^{\bar{K}_m}\subset\pi_\eta^{K_m}
= (\pi_1)_\eta^{K_m} \oplus (\pi_2)_\eta^{K_m} = (0).\end{equation*}

Suppose $m\geq c$. As in~\S\ref{sec:principal-u11}, we compute
$\dim\bar{\pi}_{\bar{\eta}}^{\bar{K}_m}$ using the fact (\ref{eq:theta})
that $\bar{\pi}_{\bar{\eta}}^{\bar{K}_m}$ is the subspace of
$\pi_\eta^{K_m}$ on which $\bar{\pi}(\theta)$ acts as the scalar
$\bar{\eta}(^s\epsilon_E^{-1})$.

Since $\pi_2 = \, ^\theta\pi_1 $ and the conjugation action of $\theta$
and $\gamma$ are the same on $G$, $\pi_1$ and $\pi_2$ form an $L$-packet
according to~\S\ref{sec:sc-sl2}. Thus $\pi = \pi_1\oplus\pi_2$ is the
restriction to $G$ of a minimal ramified supercuspidal representation
$\widetilde{\pi}$ of $\widetilde{G}$. In particular, we have an action
of $\widetilde{G}$ on $V$. Let $W$ be the one-dimensional space
$(\pi_1)_\eta^{K_c}$. Then according to the proof \cite{josh-raghuram2}
of Proposition~\ref{prop:sc-ramified-sl2}
\begin{equation*}
(\pi_1)_\eta^{K_m} = \bigoplus_{i=0}^{m-c}\widetilde{\pi}(\beta)^iW.
\end{equation*}
Now $\bar{\pi}(\theta)$ intertwines $^\theta\pi_1$ and $\pi_2$ and takes
$(\pi_1)_\eta^{K_m} = \left(^\theta\pi_1\right)_\eta^{K_m}$ to
$(\pi_2)_\eta^{K_m}$. Therefore,
\begin{equation*}
(\pi_2)_\eta^{K_m} = \bigoplus_{i=0}^{m-c} \bar{\pi} (\theta)
\widetilde{\pi}(\beta)^iW.
\end{equation*}
Let $W(i) = \widetilde{\pi}(\beta)^iW \oplus \bar{\pi}(\theta)
\widetilde{\pi}(\beta)^iW$ for $i=0, \ldots, m-c$.

Note that
\begin{align*}
\bar{\pi}(\theta)^2 = \bar{\pi}(\theta^2) &= \bar{\pi}\left( \left(\begin{array}{ll}
\epsilon_E{^s\!}\epsilon_E^{-1} &0\\
0 &\epsilon_E{^s\!}\epsilon_E^{-1} \end{array}\right) \left(\begin{array}{ll}
\epsilon_F &0\\
0 &\epsilon_F^{-1} \end{array}\right) \right)\\[.2pc]
&= \omega_{\bar{\pi}}(\epsilon_E/\,^s\epsilon_E)\bar{\pi}\left( \left(\begin{array}{ll}
\epsilon_F &0\\
0 &\epsilon_F^{-1} \end{array}\right)
\right).
\end{align*}
Thus $\bar{\pi}(\theta)^2$ acts via the scalar
$\omega_{\bar{\pi}}(\epsilon_E/\,^s\epsilon_E)\eta (\epsilon_F^{-1}) =
\bar{\eta}(^s\epsilon_E ^{-1})^2$ on $\bar{\pi}_{\eta}^{K_m}$. It
follows that $\bar{\pi}(\theta)$ exchanges the one-dimensional spaces
$\widetilde{\pi}(\beta)^iW , \bar{\pi}(\theta)\widetilde{\pi}(\beta)^iW$
since
\begin{align*}
&\bar{\pi}(\theta) (\widetilde{\pi}(\beta)^iW) = \bar{\pi}
(\theta) \widetilde{\pi}(\beta)^iW,\\[.2pc]
&\bar{\pi}(\theta) (\bar{\pi}(\theta)\widetilde{\pi}(\beta)^iW) =
\bar{\pi}(\theta)^2 (\widetilde{\pi}(\beta)^iW) =
\bar{\eta}(^s\epsilon_E ^{-1})^2\widetilde{\pi}(\beta)^iW =
\widetilde{\pi}(\beta)^iW.
\end{align*}
In particular, each $W(i)$ is stabilized by $\bar{\pi}(\theta)$.
Moreover, since $\bar{\pi}(\theta)^2$ acts via the scalar
$\bar{\eta}(^s\epsilon_E ^{-1})^2$ on $W(i)$, the eigenspaces of
$\bar{\pi}(\theta)$ on $W(i)$ corresponding to the eigenvalues
$\pm\bar{\eta}(^s\epsilon_E ^{-1})$ must each be one-dimensional. Hence
the subspace of
\begin{equation*}
\pi_\eta^{K_m} = (\pi_1)_\eta^{K_m}\oplus(\pi_2)_\eta^{K_m} =
\bigoplus_{i=0}^{m-c}W(i)
\end{equation*}
on which $\bar{\pi}(\theta)$ acts via the scalar $\bar{\eta}
(^s\epsilon_E ^{-1})$ has dimension $m-c+1$, as required.\hfill $\Box$
\end{proof}

\subsubsection*{\it \!\!\!Unramified case.}

Suppose that $(\bar{\pi},V)$ is an unramified supercuspidal
representation induced from a representation $\bar{\sigma}$ of
$\bar{K}$. It is easily seen that the restriction $\pi$ of $\bar{\pi}$
to $G$ is either
\begin{enumerate}
\renewcommand\labelenumi{(\roman{enumi})}
\leftskip .2pc
\item an irreducible unramified supercuspidal representation of $G$
induced from $K$ if the restriction of $\bar{\sigma}$ to $K$ is
irreducible, or

\item the direct sum of two irreducible unramified supercuspidal
representations of $G$ induced from $K$ if the restriction of
$\bar{\sigma}$ to $K$ is isomorphic to $\sigma_1\oplus\sigma_2$, where
$\sigma_1$ and $\sigma_2$ come from the two cuspidal representations of
$SL_2({\mathbb F}_q)$ of dimension $(q-1)/2$ (as in \S\ref{sec:sc-sl2}).
\end{enumerate}
In case (ii), we note that if $\pi$ decomposes into the direct sum of
$(\pi_1,V_1)$ and $(\pi_2, V_2)$, then $V_2 = \widetilde{\pi}
(\gamma)V_1$.

As discussed in the ramified case, if $\rho_0(\bar{\pi}) = \min\{
\rho (\bar{\pi}\otimes \chi)\}$ as $\chi$ ranges over all characters of
$\bar{G}$, then the conductors of the components of $\pi$ are
$2\rho_0(\bar{\pi}) + 2$.

\begin{theor}[(Unramified supercuspidal
representation)] \label{thm:sc-un-u11}
Let $(\bar{\pi},V)$ be an unramified supercuspidal representation of
$\bar{G}$ that is induced from $\bar{K}${\rm ,} and let
$\bar{\pi}'=\,^\alpha\bar{\pi}$. Let $\bar{\eta}$ be any character of
$\mathcal{O}_E^{\times}$ with $\bar{\eta}|_{E^1} = \omega_{\bar{\pi}}$
and $c(\bar{\eta}|_{\mathcal{O}_F^{\times}}) \leq \rho_0(\bar{\pi}) +
1$. Then $c(\bar{\pi}) = c_{\bar{\eta}}(\bar{\pi}) = 2\rho_0 (\bar{\pi})
+ 2$.
\begin{enumerate}
\renewcommand\labelenumi{\rm (\roman{enumi})}
\leftskip .2pc
\item If $\rho_0$ is odd{\rm ,} then
\begin{align*}
\hskip -1.25pc\dim (\bar{\pi}_{\bar{\eta}}^{\bar{K}_m}) &= \max\left\{
\left\lceil\frac{m-c(\bar{\pi})+1}{2}\right\rceil, 0\right\} =
\dim ((\bar{\pi}')_{\bar{\eta}}^{\bar{K}'_m}),\\[.3pc]
\hskip -1.25pc\dim (\bar{\pi}_{\bar{\eta}}^{\bar{K}'_m}) &= \max\left\{
\left\lceil\frac{m-c(\bar{\pi})-1}{2}\right\rceil, 0\right\} = \dim
((\bar{\pi}')_{\bar{\eta}}^{\bar{K}_m}).
\end{align*}

\item If $\rho_0$ is even{\rm ,} then
\begin{align*}
\hskip -1.25pc\dim ((\bar{\pi}')_{\bar{\eta}}^{\bar{K}_m}) &= \max\left\{
\left\lceil\frac{m-c(\bar{\pi})+1}{2}\right\rceil, 0\right\}= \dim
(\bar{\pi}_{\bar{\eta}}^{\bar{K}'_m}),\\[.3pc]
\hskip -1.25pc\dim ((\bar{\pi}')_{\bar{\eta}}^{\bar{K}'_m}) &= \max\left\{
\left\lceil\frac{m-c(\bar{\pi})-1}{2}\right\rceil, 0\right\} = \dim
(\bar{\pi}_{\bar{\eta}}^{\bar{K}_m}).
\end{align*}
\end{enumerate}
\end{theor}

\begin{proof}
We give a proof only for Case (ii) ($\rho_0(\bar{\pi})$ even) as the
proof for Case (i) is easily obtained therefrom by interchanging the
representations $\bar{\pi}$ and $\bar{\pi}'$. Moreover, we prove only
the first equality of each line as the second follows by conjugating by
$\alpha$.

Let $(\pi',V)$ be the restriction of $(\bar{\pi}',V)$ to $G$. Set
$c=2\rho_0(\bar{\pi}') +2$ and $\eta =
\bar{\eta}|_{\mathcal{O}_F^{\times}}$. Now $\pi'$ is a direct summand of
the restriction to $G$ of a minimal unramified supercuspidal
representation $(\widetilde{\pi},\widetilde{V})$ of $\widetilde{G}$.
Since $\widetilde{\pi}$ is unramified, it follows
from~\S\ref{sec:sc-sl2} that $^\gamma\pi'$ is isomorphic to $\pi'$ and
hence that $\widetilde{\pi}(\gamma)$ maps $V$ onto $V$. (Here we view
$V$ as a subrepresentation of $\widetilde{V}$.)

As discussed above, $\pi'$ is either an irreducible unramified
supercuspidal representation of conductor $c$ or the direct sum of two
such representations $(\pi'_1,V_1)$ and $(\pi'_2,V_2)$, where $V_2 =
\widetilde{\pi}(\gamma)V_1$. By
Propositions~\ref{prop:sc-unramified-2-sl2} and
\ref{prop:sc-unramified-4-sl2}, the level $l$ of the inducing data for
these representations is $c/2 = \rho_0(\bar{\pi}') + 1$. As in the
ramified case, we have
$\dim\left(\bar{\pi}'\right)_{\bar{\eta}}^{\bar{K}_m}=0$ if $m< c$.

Suppose $m\geq c$. By (\ref{eq:theta}), to find
$\dim\left(\bar{\pi}'\right)_{\bar{\eta}}^{\bar{K}_m}$, we compute the
dimension of the subspace of $(\pi')_\eta^{K_m}$ on which
$\bar{\pi}'(\theta)$ acts as the scalar $\bar{\eta}(^s\epsilon_E^{-1})$.
Let $W=(\pi')_\eta^{K_c}$. Since $l=\rho_0(\bar{\pi}') + 1$ is odd, it
follows from Propositions~\ref{prop:sc-unramified-2-sl2} and
\ref{prop:sc-unramified-4-sl2} and their proofs~\cite{josh-raghuram2}
that $\dim (W)=2$ and
\begin{equation*}
(\pi')_\eta^{K_m} = \bigoplus_{i = 0}^{\lfloor (m-c)/2\rfloor}\pi \left(
\left(\begin{array}{ll}
\varpi_F &0\\
0 &\varpi_F^{-1} \end{array}\right) \right)^iW.
\end{equation*}

In fact, from the proof of Proposition 3.3.4 in~\cite{josh-raghuram2},
it follows that for a certain vector $\phi \in W$, $W = \mathbb{C}\phi
\oplus \mathbb{C}\widetilde{\pi}(\gamma)\phi$. If $\pi'$ is irreducible,
then since conjugation by $\gamma$ and $\theta$ have the same effect on
$G$, $\widetilde{\pi}(\gamma)$ and $\bar{\pi}'(\theta)$ are both
elements of the one-dimensional space $\text{Hom}(^\theta\pi' ,\pi')$.
They are therefore equal up to scalars so $W = \mathbb{C}\phi \oplus
\mathbb{C}\bar{\pi}'(\theta)\phi$. If $\pi'$ is reducible, then we may
further assume that $\phi \in W\cap V_1$ by
Proposition~\ref{prop:sc-unramified-4-sl2-test}. In this case,
$\widetilde{\pi}(\gamma)$ and $\bar{\pi}'(\theta)$ are both elements of
the one-dimensional space $\text{Hom}(^\theta\pi'_1 ,\pi'_2)$ so $W =
\mathbb{C}\phi \oplus \mathbb{C}\bar{\pi}'(\theta)\phi$ as above.

As in the ramified case, $\bar{\pi}'(\theta)^2$ acts via the scalar
$\bar{\eta}(^s\epsilon_E ^{-1})^2$ on
$\left(\bar{\pi}'\right)_{\eta}^{K_m}$. It follows that
$\bar{\pi}'(\theta)$ exchanges the one-dimensional spaces
$\mathbb{C}\phi , \bar{\pi}'(\theta)\mathbb{C}\phi$ since
\begin{align*}
&\bar{\pi}'(\theta)(\mathbb{C}\phi) = \mathbb{C}\bar{\pi}'(\theta) \phi,\\[.2pc]
&\bar{\pi}'(\theta)\left(\mathbb{C}\bar{\pi}'(\theta) \phi\right) =
\bar{\pi}'(\theta)^2\left(\mathbb{C}\phi\right) =
\bar{\eta}(^s\epsilon_E ^{-1})^2\mathbb{C}\phi =
\mathbb{C}\phi.
\end{align*}
In particular, $\bar{\pi}'(\theta)$ stabilizes $W$. Again as in the
ramified case, these facts imply that the eigenspaces of
$\bar{\pi}'(\theta)$ on $W$ corresponding to the eigenvalues
$\pm\bar{\eta}(^s\epsilon_E ^{-1})$ must each be one-dimensional. The
same is clearly true of $\pi \left( \left(\begin{smallmatrix} \varpi_F
&0\\ 0 &\varpi_F^{-1} \end{smallmatrix}\right) \right)^i W$ for $i =
0,\ldots ,\lfloor (m-c)/2\rfloor$. It follows that the subspace of
$(\pi')_\eta^{K_m}$ on which $\bar{\pi}'(\theta)$ acts via the scalar
$\bar{\eta}(^s\epsilon_E ^{-1})$ has dimension $\lceil (m-c+1)/2\rceil$
as required.

The computation of $\dim \left( \bar{\pi}'
\right)_{\bar{\eta}}^{\bar{K}'_m}$ is entirely analogous.\hfill $\Box$
\end{proof}

Now suppose that $\bar{\pi}$ is a supercuspidal representation of
$\bar{G}$ of conductor $c$. We consider the effect of a Whittaker
functional $\Lambda_\psi$ on $\bar{\pi}_{\bar{\eta}}^{\bar{K}_c},
\bar{\pi}_{\bar{\eta}}^{\bar{K}'_c}$. For this we need to choose the
character $\bar{\eta}$ somewhat carefully. Let $\bar{\Pi}$ be the
$L$-packet of $\bar{G}$ containing $\bar{\pi}.$ Then the restriction to
$G$ of the direct sum of representations in $\bar{\Pi}$ is also the
restriction to $G$ of a minimal supercuspidal representation
$\widetilde{\pi}$ coming via Kutzko's construction. We require
$\bar{\eta}|_{\mathcal{O}_F^{\times}} = \omega_{\widetilde{\pi}}.$

\begin{propo} {\rm (Test vectors for supercuspidal representations)}
$\left.\right.$\vspace{.5pc}

\noindent Suppose that $\bar{\pi}$ is an irreducible supercuspidal
representation of $\bar{G}$ of conductor $c$. Let $\bar{\eta}$ be a
character of $\mathcal{O}_E^{\times}$ with
$\bar{\eta}|_{E^1}=\omega_{\bar{\pi}}$ and
$\bar{\eta}|_{\mathcal{O}_F^{\times}} = \omega_{\widetilde{\pi}}$
{\rm (}see above{\rm )}. Let $\psi=\psi_F$.
\begin{enumerate}
\renewcommand\labelenumi{\rm (\roman{enumi})}
\leftskip .2pc
\item If $\bar{\pi}$ is ramified{\rm ,} then $\bar{\pi}$ is $\psi$-generic.
Moreover{\rm ,} $\Lambda_\psi (v)\neq 0$ for all non-zero $v$ in
$\bar{\pi}_{\bar{\eta}}^{\bar{K}_c}$ or
$\bar{\pi}_{\bar{\eta}}^{\bar{K}'_c}$.

\item If $\bar{\pi}$ is unramified and induced from $\bar{K}${\rm ,} let
$\bar{\pi}' = \,^\alpha\bar{\pi}$.

%\begin{enumerate}
%\renewcommand\labelenumii{\rm (\alph{enumii})}
\leftskip .2pc
(a) If $\rho_0 (\bar{\pi}) = \rho_0 (\bar{\pi}')$
is odd{\rm ,} then $\bar{\pi}$ is $\psi$-generic and $\bar{\pi}'$
is $\psi_\varpi$-generic. Moreover{\rm ,} $\Lambda_\psi (v)\neq 0$
for all non-zero $v\in \bar{\pi}_{\bar{\eta}}^{\bar{K}_c}$ and
$\Lambda_{\psi_\varpi} (v)\neq 0$ for all non-zero\break $v\in
\left(\bar{\pi}'\right)_{\bar{\eta}}^{\bar{K}'_c}$.

(b) If $\rho_0 (\bar{\pi}) = \rho_0 (\bar{\pi}')$ is even{\rm ,}
then $\bar{\pi}'$ is $\psi$-generic and $\bar{\pi}$ is
$\psi_\varpi$-generic. Moreover{\rm ,} $\Lambda_\psi (v)\neq 0$
for all non-zero $v\in
\left(\bar{\pi}'\right)_{\bar{\eta}}^{\bar{K}_c}$ and
$\Lambda_{\psi_\varpi} (v)\neq 0$ for all non-zero $v\in
\bar{\pi}_{\bar{\eta}}^{\bar{K}'_c}$.
%\end{enumerate}
\end{enumerate}\vspace{-1pc}
\end{propo}

\begin{proof}
Let $\pi$ be the restriction of $\bar{\pi}$ to $G$. As in
Proposition~\ref{prop:ps-test-u11}, we note that the restriction of
$\Lambda_\psi$ to any $\psi$-generic component of $\pi$ is a
$\psi$-Whittaker functional on that component, while its restriction to
any non-$\psi$-generic component is $0$.

Suppose first that $\bar{\pi}$ is ramified (case (i)). Then $\pi$
decomposes as the direct sum $\pi_1\oplus\pi_2$ of irreducible ramified
supercuspidal representations of conductor $c$. By
Proposition~\ref{prop:sc-ramified-sl2-test}, only one summand, say
$\pi_1$, is $\psi$-generic and we have that $\Lambda_{\psi}$ is non-zero
on $(\pi_1)^{K_c}_{\eta}.$ Now $\bar{\pi}_{\bar{\eta}}^{\bar{K}'_c}$ is
the space of vectors in $\pi_\eta^{K_c} =
\left(\pi_1\right)_\eta^{K_c}\oplus\left(\pi_2\right)_\eta^{K_c}$ on
which $\bar{\pi}(\theta)$ acts as the scalar
$\bar{\eta}(^s\epsilon_E^{-1})$. As observed in the proof of
Theorem~\ref{thm:sc-ram-u11}, $\bar{\pi}(\theta)$ exchanges
$\left(\pi_1\right)_\eta^{K_c} $ and $\left(\pi_2\right)_\eta^{K_c}$.
Therefore, $\bar{\pi}_{\bar{\eta}}^{\bar{K}_c}$ cannot lie in either
$\left(\pi_1\right)_\eta^{K_c}$ or $\left(\pi_2\right)_\eta^{K_c}$. In
particular, if $v \in \bar{\pi}_{\bar{\eta}}^{\bar{K}_c}$ is written as
$v_1+v_2$ with $v_i\in\left(\pi_i\right)_\eta^{K_c}$, then $v_1,v_2\neq
0$. Since $\pi_1 $ is $\psi$-generic and $\pi_2$ is not, we get
\begin{equation*}
\Lambda_\psi (v) = \Lambda_\psi (v_1) \neq 0.
\end{equation*}

We now give a proof in case (ii). We only prove (a) as the proof of (b)
is obtained by interchanging $\bar{\pi}$ and $\bar{\pi}'$.

Suppose that $\bar{\pi}$ is unramified and induced from $\bar{K}$ and
that $\rho_0 (\bar{\pi})$ is odd. Then $\pi$ is also unramified, induced
from $K$, and has conductor $c$. As noted in the proof of
Theorem~\ref{thm:sc-un-u11}, $\pi_\eta^{K_c} =
\mathbb{C}\phi_1\oplus\mathbb{C}\bar{\pi}(\theta)\phi_1.$ Since the
level of the inducing data of $\pi$ is $c/2$, which is even, $\pi$ is
$\psi$-generic by Proposition~\ref{prop:sc-unramified-2-sl2-test}.
Moreover, $\Lambda_\psi (\phi_1)\neq 0$, while $\Lambda_\psi
(\bar{\pi}(\theta)\phi_1) = \Lambda_\psi (\phi_\epsilon ) = 0$. By the
proof of Theorem~\ref{thm:sc-ram-u11}, $\bar{\pi}(\theta)$ exchanges
$\mathbb{C}\phi_1$ and $\mathbb{C}\bar{\pi}(\theta)\phi_1$. Therefore,
just as in the ramified case, if $v = a\phi_1+b\bar{\pi}(\theta)\phi_1$,
then $a,b\neq 0$. It follows that
\begin{equation*}
\Lambda_\psi (v) = a\Lambda_\psi (\phi_1) \neq 0.
\end{equation*}
The proof of the non-vanishing of $\Lambda_{\psi_{\varpi}}$ is entirely
analogous.\hfill $\Box$
\end{proof}

\begin{rem}{\rm
We have only considered the unitary group $U(1,1)$ for an unramified
extension $E/F.$ The entire series of results go through with some minor
modifications if instead we considered ramified extensions.}
\end{rem}

\subsection{\it Comparison of conductor with other invariants}
\label{sec:comparison-u11}

\setcounter{theore}{0}
\begin{theor}[(Relation of conductor with other invariants for
$\bar{\hbox{\bfseries\itshape G}}$)] \label{thm:conductor-otherinvariants-u11}
Let $\bar{\pi}$ be an irreducible admissible supercuspidal
representation of $\bar{G}.$ The relation between its conductor
$c(\bar{\pi})$ and its minimal depth $\rho_0(\bar{\pi})$ is given
by
\begin{equation*}
\rho_0(\bar{\pi}) = \frac{c(\bar{\pi})-2}{2}.
\end{equation*}
If $\pi$ is an irreducible subrepresentation of the restriction of
$\bar{\pi}$ to $G$ then
\begin{equation*}
c(\bar{\pi}) = c(\pi).
\end{equation*}
\end{theor}

\begin{proof}
This follows from Theorems \ref{thm:sc-ram-u11} and
\ref{thm:sc-un-u11}.\hfill $\Box$
\end{proof}

\section{Towards multiplicity one for newforms}
\label{sec:multiplicity-one}

Given an irreducible representation $\bar{\pi}$ of $\bar{G}$ and a
character $\bar{\eta}$ of $\mathcal{O}_E^{\times}$ such that
$c_{\bar{\eta}}(\bar{\pi}) = c(\bar{\pi})$, one can ask if we have $\dim
(V^{\bar{K}_{c(\bar{\pi})}}_{\bar{\eta}}) = 1$. The answer is that this
is often the case but is not true in general. Indeed, we have $\dim
(V^{\bar{K}_{c(\bar{\pi})}}_{\bar{\eta}}) = 1$ unless $\bar{\pi}$ is the
principal series representation $\bar{\pi}(\bar{\chi})$, where
$\bar{\chi}$ is ramified and $\bar{\chi}|_{\mathcal{O}_E^{\times}} =
\,^s\bar{\chi}^{-1}|_{\mathcal{O}_E^{\times}}$. For these exceptional
representations, the dimension of the space of newforms is two.

Nevertheless, in all cases we have proved that an appropriate Whittaker
functional is non-vanishing on some newform. This can be used to
formulate a kind of a multiplicity one result if we consider the
quotient of the space  of newforms %$V^{\bar{K}_{c(\bar{\pi})}}_{\bar{\eta}}$
by the kernel of this Whittaker functional. More precisely, if $\Psi$ is
a non-trivial additive character of $F$ of conductor either
$\mathcal{O}_F$ or $\mathcal{P}_F^{-1}$ such that $\bar{\pi}$ is
$\Psi$-generic, and $\Lambda_{\Psi}$ is a $\Psi$-Whittaker functional
for $\bar{\pi}$, then we have
\begin{equation*}
\dim \left( \frac{V^{\bar{K}_{c(\bar{\pi})}}_{\bar{\eta}}}
{V^{\bar{K}_{c(\bar{\pi})}}_{\bar{\eta}} \cap
{\rm kernel}(\Lambda_{\Psi})} \right) = 1.
\end{equation*}

Another possibility is to consider some canonical non-degenerate
bilinear form on the space $V^{\bar{K}_{c(\bar{\pi})}}$ and consider the
orthogonal complement of the subspace
$V^{\bar{K}_{c(\bar{\pi})}}_{\bar{\eta}} \cap {\rm
kernel}(\Lambda_{\Psi})$ as a candidate for a one-dimensional space of
newforms. Then the multiplicity one result is formulated as
\begin{equation*}
\dim \left(V^{\bar{K}_{c(\bar{\pi})}}_{\bar{\eta}} \cap
{\rm kernel}(\Lambda_{\Psi})^{\perp}\right) = 1.
\end{equation*}

\section*{Acknowledgments}

We thank Benedict Gross, Dipendra Prasad and Paul~Sally~Jr. for some
helpful correspondence. We thank Brooks Roberts for his interest in this
work and for his generous comments on new avenues of thought that might
stem from it. The second author is grateful for the warm and pleasant
working experience at Bucknell University where a substantial part of
this work was completed. We thank the referee for a very careful reading
of the manuscript.


\begin{thebibliography}{99}
\bibitem{casselman} Casselman~W, On some results of Atkin and Lehner,
{\it Math. Ann.} {\bf 201} (1973) 301--314

\bibitem{gelbart-knapp} Gelbart~S and Knapp~A~W,
$L$-indistinguishability and $R$ groups for the special linear group,
{\it Adv. Math.} {\bf 43(2)} (1982) 101--121

\bibitem{gross-prasad} Gross~B and Prasad~D, Test vectors for linear
forms, {\it Math. Ann.} {\bf 291(2)} (1991) 343--355

\bibitem{jacquet-ps-shalika} Jacquet~H, Piatetski-Shapiro~I and
Shalika~J, Conducteur des representations du groupe linaire, {\it Math.
Ann.} {\bf 256(2)} (1981) 199--214

\bibitem{kutzko2} Kutzko~P~C, On the supercuspidal representations of
$Gl\sb{2}$, {\it Am. J. Math.} {\bf 100(1)} (1978) 43--60

\bibitem{kutzko3} Kutzko~P~C, On the supercuspidal representations of
$Gl\sb{2}$. II, {\it Am. J. Math.} {\bf 100(4)} (1978) 705--716

\bibitem{kutzko4} Kutzko~P~C, The exceptional representations of
$Gl\sb{2}$, {\it Compositio Math.} {\bf 51(1)} (1984) 3--14

\bibitem{kutzko-sally} Kutzko~P~C and Sally Jr.~P, All supercuspidal
representations of $SL\sb{l}$ over a $p$-adic field are induced,
Representation theory of reductive groups (Utah: Park City) (1982)
pp.~185--196; {\it Prog. Math.} (Boston, MA: Birkha\"user)
(1983) vol.~40

\bibitem{labesse-langlands} Labesse~J-P and Langlands~R~P,
$L$-indistinguishability for $SL(2)$, {\it Canad. J. Math.} {\bf 31(4)}
(1979) 726--785

\bibitem{josh-raghuram1} Lansky~J and Raghuram~A, A remark on the
correspondence of representations between $GL(n)$ and division algebras,
{\it Proc. Am. Math. Soc.} {\bf 131(5)} (2003) 1641--1648

\bibitem{josh-raghuram2} Lansky~J and Raghuram~A, Conductors and
newforms for $SL(2)$, (submitted) Preprint available at
http://www.math.uiowa.edu/\~{}araghura/newforms\_sl2.dvi

\bibitem{12} Mann~W~R, Local level raising for $GL(n)$, Ph.D. thesis
(Harvard University) (2001)

\bibitem{moy-prasad1} Moy~A and Prasad~G, Unrefined minimal $K$-types
for $p$-adic groups, {\it Invent. Math.} {\bf 116(1--3)} (1994) 393--408

\bibitem{moy-sally} Moy~A and Sally Jr~P, Supercuspidal representations
of $SL\sb{n}$ over a $p$-adic field: the tame case, {\it Duke Math. J.}
{\bf 51(1)} (1984) 149--161

\bibitem{dipendra} Prasad~D, Trilinear forms for representations of
$GL(2)$ and local $\epsilon$-factors, {\it Compositio Math.} {\bf 75(1)}
(1990) 1--46

\bibitem{dipendra-raghuram} Prasad~D and Raghuram~A, Kirillov theory for
$GL\sb 2(\mathcal{D})$ where $\mathcal{D}$ is a division algebra over a
non-Archimedean local field, {\it Duke Math. J.} {\bf 104(1)} (2000)
19--44

\bibitem{rogawski} Rogawski~J, Automorphic representations of unitary
groups in three variables, {\it Ann. Math. Stud.} 123 (Princeton
University Press) (1990)

%\bibitem{russ} Russell Mann~W, Local level raising for $GL(n)$, Harvard
%Thesis (2001)

\bibitem{schmidt} Schmidt~R, Some remarks on local newforms for $GL(2)$,
{\it J. Ramanujan Math. Soc.} {\bf 17(2)} (2002) 115--147

\bibitem{shimizu} Shimizu~H, Some examples of new forms, {\it J. Fac.
Sci. Univ. Tokyo Sect. IA Math.} {\bf 24(1)} (1977) 97--113
\end{thebibliography}
\end{document}